\documentclass[a4paper, 11pt]{article}

\usepackage{longtable}
\usepackage{geometry}
\usepackage{graphicx}
\usepackage{subfig}
\usepackage{color}
\usepackage{xcolor}
\usepackage{bm}
\usepackage{fixmath}
\usepackage{mathrsfs}
\usepackage{ulem}
\usepackage{multirow}
\usepackage{pdflscape}
\usepackage{rotating}
\usepackage{lscape}
\usepackage{pdflscape}
\usepackage{titlesec}
\usepackage{float}
\usepackage{hyperref}
\usepackage{tabularx}
\usepackage{mathabx}
\usepackage{wasysym}

\titleclass{\subsubsubsection}{straight}[\subsection]

\newcounter{subsubsubsection}[subsubsection]
\renewcommand\thesubsubsubsection{\thesubsubsection.\arabic{subsubsubsection}}

\titleformat{\subsubsubsection}
  {\normalfont\normalsize\bfseries}{\thesubsubsubsection}{1em}{}
\titlespacing*{\subsubsubsection}
{0pt}{3.25ex plus 1ex minus .2ex}{1.5ex plus .2ex}

\makeatletter

\renewcommand\paragraph{\@startsection{paragraph}{5}{\z@}%
  {3.25ex \@plus1ex \@minus.2ex}%
  {-1em}%
  {\normalfont\normalsize\bfseries}}

\renewcommand\subparagraph{\@startsection{subparagraph}{6}{\parindent}%
  {3.25ex \@plus1ex \@minus .2ex}%
  {-1em}%
  {\normalfont\normalsize\bfseries}}

\def\toclevel@subsubsubsection{4}
\def\toclevel@paragraph{5}
\def\toclevel@paragraph{6}
\def\l@subsubsubsection{\@dottedtocline{4}{7em}{4em}}
\def\l@paragraph{\@dottedtocline{5}{10em}{5em}}
\def\l@subparagraph{\@dottedtocline{6}{14em}{6em}}

\makeatother

\def\b{\begin{eqnarray}}
\def\e{\end{eqnarray}}

\def\n{\noindent}

\def\1{\vskip.1cm}
\def\2{\vskip.2cm}
\def\3{\vskip.3cm}
\def\5{\vskip.5cm}

\makeatletter
\newcommand*\bigdot{\mathpalette\bigdot@{.5}}
\newcommand*\bigdot@[2]{\mathbin{\vcenter{\hbox{\scalebox{#2}{$\m@th#1\bullet$}}}}}
\makeatother

\begin{document}

\begin{center}
{\huge \textbf{A Method for Locating the Real Roots  \vskip.1cm of the Symbolic Quintic Equation Using \vskip.25cm Quadratic Equations}}

\vspace{9mm}
\noindent
{\large \bf Emil M. Prodanov} \vskip.4cm
{\it School of Mathematical Sciences, Technological University Dublin,
\vskip.1cm
Park House, Grangegorman, 191 North Circular Road, Dublin
D07 EWV4, Ireland,}
\vskip.1cm
{\it e-mail: emil.prodanov@tudublin.ie} \\
\vskip.5cm
\end{center}

\vskip2cm

\begin{abstract}
\n
A method is proposed with which the locations of the roots of the monic symbolic quintic polynomial $x^5 + a_4 x^4 + a_3 x^3 + a_2 x^2 + a_1 x + a_0$ can be determined using the roots of two {\it resolvent} quadratic polynomials: $q_1(x) = x^2 + a_4 x + a_3$ and $q_2(x) = a_2 x^2 + a_1 x + a_0$, whose coefficients are exactly those of the quintic polynomial. The different cases depend on the coefficients of $q_1(x)$ and $q_2(x)$ and on some specific relationships between them. The method is illustrated with the full analysis of one of the possible cases. Some of the roots of the symbolic quintic equation for this case have their isolation intervals determined and, as this cannot be done for all roots with the help of quadratic equations only, finite intervals containing 1 or 3 roots, or 0 or 2 roots, or, rarely, 0, or 2, or 4 roots of the quintic are identified. Knowing the stationary points of the quintic polynomial, lifts the latter indeterminacy and allows one to find the isolation interval of each of the roots of the quintic. Separately, using the complete root classification of the quintic, one can also lift this indeterminacy. The method also allows to see how variation of the coefficients of the quintic affect its roots. No root finding iterations or any numerical approximations are used and no equations of degree higher than 2 are solved.
\end{abstract}


\newpage

\subparagraph{\hskip-.6cm 1 \hskip0.2cm Introduction  \vskip0.5cm}
\hskip-1cm As it is neither possible to get the roots of an equation of degree five or higher in terms of radicals (Abel--Ruffini theorem) nor it is possible to find a closed-form solution of such equations [1], the localization of the real roots of the quintic equation by determination of their isolation intervals or by finding finite intervals containing $m - 2n \ge 1$ real roots (where $m \le 4$ and $n = 0, 1, 2$) would prove to be a valuable benefit when one deals with a quintic equation. This is especially important in situations in which the coefficients of the quintic depend on the parameters of the phenomenon modeled by the equation --- as numerical methods become rather involved in such cases. \\
For example, in celestial mechanics, in the closed Sun--Earth system, for a test body at the Lagrange point $L_1$ (which is along the axis Sun--Earth, about $1.5 \times 10^6$ km from Earth's centre inside Earth's orbit), the gravitational field of the Earth balances that of the Sun. The orbital period of the test body will be exactly equal to the orbital period of the Earth. Hence, objects at $L_1$ tend to ``stay put" (the libration point $L_1$ is unstable). The point $L_1$ provides an uninterrupted view of the Sun and the Solar and Heliospheric Observatory Satellite (SOHO) has been at the $L_1$ point since 1995. The distance between the centre of the Earth and the Lagrange point $L_1$ is the only positive root of a quintic equation whose coefficients involve the reduced mass of the Earth (in units of the total mass of the system) and the assumed constant distance between the Earth and the Sun (there are four other Lagrange points in the Sun--Earth system). \\
Root classification and complete root classification of parametric polynomials have been extensively studied --- see [2] and the references therein. The root classification provides the collection of all possible cases of the polynomial roots and consists of a list of the multiplicities of all roots (real and complex). The complete root classification consists of the root classification, together with the conditions which the equation coefficients should satisfy for each of the cases of the root classification. Neither the root classification, nor the complete root classification deal with the location of the roots of the polynomial. \\
The complete root classification for the depressed monic quintic, $x^5 + p x^3 + q x^2 + r x + s$, is proposed in 1996 by Yang, Xiaorong, and Zhen [2]. It is as follows. Introducing:
\b 
\label{de2}
D_2 \!\!\! & = & \!\!\! -p, \\
D_3 \!\!\! & = & \!\!\!  40 r p  - 12 p^3 - 45q^2, \\
D_4 \!\!\! & = & \!\!\! 12 p^4 r - 4 p^3 q^2 + 117 p r q^2 - 88 r^2 p^2 - 40 p^2 q s + 125 p s^2 - 27 q^4 - 300 q r s + 160 r^3, \nonumber \\ \\
D_5 \!\!\! & = & \!\!\! - 1600 q s r^3 - 3750 p s^3 q + 2000 p s^2 r^2 - 4 p^3 q^2 r^2 + 16 p^r q^3 s - 900 r s^2 p^3 \nonumber \\
& & + \, 825 p^2 q^2 s^2 + 144 p q^2 r^3 + 2250 q^2 r s^2 + 16 r^4 p^3 + 108 p^5 s^2 - 128 r^4 p^2 - 27 q^4 r^2 \nonumber \\
& & + \, 108 q^5 s + 256 r^3 + 3125 s^4 - 72 p^4 r s q + 560 p^2 r^2 s q - 630 p r s q^4, \\
E_2 \!\!\! & = & \!\!\! 160 r^2 p^3 + 900 q^2 r^2 - 48 r p^5 + 60 q^2 p^2 r + 1500 p q r s + 16 q^2 p^4 - 1100 q p^3 s \nonumber \\
& & + \, 625 s^2 p^2 - 3375 q^3 s, \\
\label{ef2}
F_2 \!\!\! & = & \!\!\! 3 q^2 - 8 r p,
\e 
one has:
\begin{center}
\begin{tabular}{rlr}
1 & $D_5 > 0 \wedge D_4 > 0 \wedge D_3 > 0 \wedge D_2 > 0$ & $\{1, 1, 1, 1, 1 \}$ \\
2 & $D_5 > 0 \wedge (D_4 \le 0 \vee D_3 \le 0 \vee D_2 \le 0)$ & $\{1 \}$ \\
3 & $D_5 < 0$ & $\{1, 1, 1 \}$ \\
4 & $D_5 = 0 \wedge D_4 > 0$ & $\{2, 1, 1, 1 \}$ \\
5 & $D_5 = 0 \wedge D_4 < 0$ & $\{2, 1 \}$ \\
6 & $D_5 = 0 \wedge D_4 = 0 \wedge D_3 > 0 \wedge E_2 \ne 0$ & $\{2, 2, 1 \}$ \\
7 & $D_5 = 0 \wedge D_4 = 0 \wedge D_3 > 0 \wedge E_2 = 0$ & $\{3, 1, 1 \}$ \\
8 & $D_5 = 0 \wedge D_4 = 0 \wedge D_3 < 0 \wedge E_2 \ne 0$ & $\{1 \}$ \\
9 & $D_5 = 0 \wedge D_4 = 0 \wedge D_3 < 0 \wedge E_2 = 0$ & $\{3 \}$ \\
10 & $D_5 = 0 \wedge D_4 = 0 \wedge D_3 = 0 \wedge D_2 \ne 0 \wedge F_2 \ne 0 $ & $\{3, 2 \}$ \\
11 & $D_5 = 0 \wedge D_4 = 0 \wedge D_3 = 0 \wedge D_2 \ne 0 \wedge F_2 = 0 $ & $\{4, 1 \}$ \\
12 & $D_5 = 0 \wedge D_4 = 0 \wedge D_3 = 0 \wedge D_2 = 0 $ & $\{5 \}$ \\
\end{tabular}
\end{center}
The lists in the figure brackets in the third column, in Yang, Xiaorong, and Zhen notation, show the multiplicities of the real roots. For example, $\{2, 1, 1, 1\}$ means five real roots: one double and three simple. \\
The goal of this work is to provide a tool for finding the {\it isolation intervals} of some of the real roots of the monic symbolic quintic equation $x^5 + a_4 x^4 + a_3 x^3 + a_2 x^2 + a_1 x + a_0 = 0$ and, for the roots for which this is not possible, to allow the determination of finite intervals containing clusters of roots. Despite the complexity of the quintic equation, this can be achieved with the help of the roots of two {\it resolvent quadratics}, $q_1(x) = x^2 + a_4 x + a_3$ and $q_2(x) = - a_2 x^2 - a_1 x - a_0$, and the end-point of the isolation or clustering intervals will turn out to be the roots of these two quadratics (some other resolvent quadratic equations will also be used in the analysis). The coefficients of the resolvent quadratics are exactly those of the original quintic. The parabola $q_2(x)$ is viewed as an element of a congruence of parabolas $- a_2 x^2 - a_1 x - \alpha$, which differ from each other by their free term and which foliate the $xy$-plane. The localization of the roots is achieved by ``splitting" the quintic equation into an equation for the intersection points of $x^3 q_1(x)$ and $q_2(x)$.  The analysis is very easy due to the fact that $x = 0$ is (at least) a triple root of the quintic $x^3 q_1(x)$ and this makes the curve $x^3 q_1(x)$ very easy to study, alongside the congruence of parabolas. The localization of the roots of the quintic equation also allows to see how the coefficients of the equation affect the roots. \\
As can be seen from the above complete root classification, despite the complexity of the conditions for the various cases, one can eliminate all of the indeterminacy associated with the method proposed in this work by depressing the quintic ($x \to x - a_4/5$) and studying the resulting $D_2, \,\, D_3, \,\, D_4, \,\, D_5, \,\, E_2,$ and $F_2$. Hence, one will know the {\it exact} number of roots within the intervals determined by the roots of the two {\it resolvent} quadratic polynomials: $q_1(x) = x^2 + a_4 x + a_3$ and $q_2(x) = a_2 x^2 + a_1 x + a_0$. \\
No quartic or cubic equations will be solved in this work, even though finding the stationary points of the given quintic will allow the determination of the isolation interval of each root of the quintic. Also, no root finding iterations or numerical methods of any kind will be used. \\

\newpage
\subparagraph{\hskip-.6cm 2 \hskip0.2cm The Method \vskip0.5cm}
\hskip-1cm Consider the general monic quintic equation
\b
\label{quintic}
Q(x) \equiv x^5 + a_4 x^4 + a_3 x^3 + a_2 x^2 + a_1 x + a_0 = 0
\e
whose coefficients $a_0, \, a_1, \, a_2, a_3,$ and $a_4$ are fixed given real numbers (whatever values they may have) and whose real roots are denoted by $x_i$ (these can be 1, or 3, or 5). Let $\xi_i$ denote the stationary points of the quintic (these can be 0, or 2, or 4). There is just one real root of the quintic if it has no stationary points, or there are one or three real roots in case of two stationary points, or there are one, or there, or five real roots (counted with their multiplicities) in case of four stationary points. \\
A key element of the method is to ``release" the fixed given coefficient $a_0$ so that it can be varied. In this way, the given quintic $x^5 + a_4 x^4 + a_3 x^3 + a_2 x^2 + a_1 x + a_0$ (with fixed $a_0$, along all other coefficients) can be viewed as an element of a one-parameter congruence of quintics which foliate the $xy$-plane. For the curves of this congruence,  $a_1, \, a_2, a_3,$ and $a_4$ are all fixed and the different quintics differ from each other by their free term only. Within the congruence, there are several privileged quintics. The one for which the free term is zero, will be referred to as {\it separatrix quintic}, namely, this is the quintic $Q_0(x) = x^5 + a_4 x^4 + a_3 x^3 + a_2 x^2 + a_1 x$ which passes through the origin. In addition to the separatrix quintic, there are further 0, or 2, or 4 privileged quintics $Q_i(x)$ (counted with their multiplicities). These are the quintics for which the abscissa is tangent to their graphs at the stationary point $\xi_i$. Note that all quintics from the congruence have the same set of stationary points. As will be addressed further, the free term of each of the privileged quintics $Q_i(x)$, with $i > 0,$ is given by $\alpha_i = a_0 - Q(\xi_i)$ --- see Figure 1a. \\
Each of the privileged quintics $Q_i(x)$ has (at least) a double root at $\xi_i$. That is, $Q_i(\xi_i) = 0$ and $Q_i'(\xi_i) = 0$ both hold. \\
Another key element of the method, based on the analysis of [3], is to view the given quintic equation (\ref{quintic}) as
\b
\label{split}
x^3 q_1(x) = q_2(x),
\e
where $x^3 q_1(x) \equiv x^3 (x^2 + a_4 x + a_3)$ and $q_2(x) \equiv  - a_2 x^2 - a_1 x - a_0$. Then the real roots $x_i$ of the quintic equation (\ref{quintic}) are the intersection points of these two curves. \\
The two quadratic equations
\b
q_1(x) & \equiv & x^2 + a_4 x + a_3 = 0, \\
- q_2(x) & \equiv & a_2 x^2 + a_1 x + a_0 = 0,
\e
will be referred to as {\it first and second resolvent quadratic equations}. The roots of these quadratic equations will be used for the localization of the roots of the quintic equation. \\
Note that in the quadratic $q_2(x)$, the coefficient $a_0$ can still vary. In this manner, the $xy$-plane is foliated by a congruence of parabolas with different free term (or straight lines, should $a_2$ happen to be zero). The curve $x^3 q_1(x)$ intersects the different parabolas of this  congruence at either 1, or 3, or 5 points. As each privileged quintic $Q_i(x)$ has (at least) a double root, ``splitting" the privileged quintics into two components, $x^3 q_1(x)$ and $- a_2 x^2 - a_1 x - \alpha_i$, and equating them, leads to up to four quintic equations each of which has a root of order (at least) two. These double roots are exactly the stationary points $\xi_i$ of the original quintic equation and at these stationary points, both $x^3 q_1(x) = - a_2 x^2 - a_1 x - \alpha_i$ and $[x^3 q_1(x)]' = (- a_2 x^2 - a_1 x - \alpha_i)'$ hold. Namely, the curves $x^3 q_1(x)$ and  $- a_2 x^2 - a_1 x - \alpha_i$ have the same value at the stationary point $\xi_i$ and, further, at the stationary point $\xi_i$, the tangents to the two curves coincide --- see Figure 1b. \\
The second of the above two simultaneous equations,
\b
\label{quartic}
x^4 + a x^3 + b x^2 + c x + d = 0,
\e
with $a = (4/5) a_4$, $b = (3/5) a_3$, $c = (2/5) a_2$, and $d = (1/5) a_1$ is the equation for the stationary points $\xi_i$ of the original quintic. This equation will be referred to as {\it auxiliary quartic equation}. It is a general quartic equation and the locations of its roots $\xi_i$ can be studied with the analysis presented in [4] --- with the help of its own set of corresponding {\it resolvent} and {\it auxiliary equations}. \\
If one solves explicitly the quartic equation (\ref{quartic}) to determine the stationary points $\xi_i$ of the quintic, then the privileged quintics $Q_i(x)$ can be immediately identified by the determination of their free terms with the former of the two simultaneous equations:
\b
\alpha_i = a_0 - Q_i(\xi_i),
\e
(as mentioned earlier). 
\begin{center}
\begin{tabular}{cc}
\includegraphics[width=67mm]{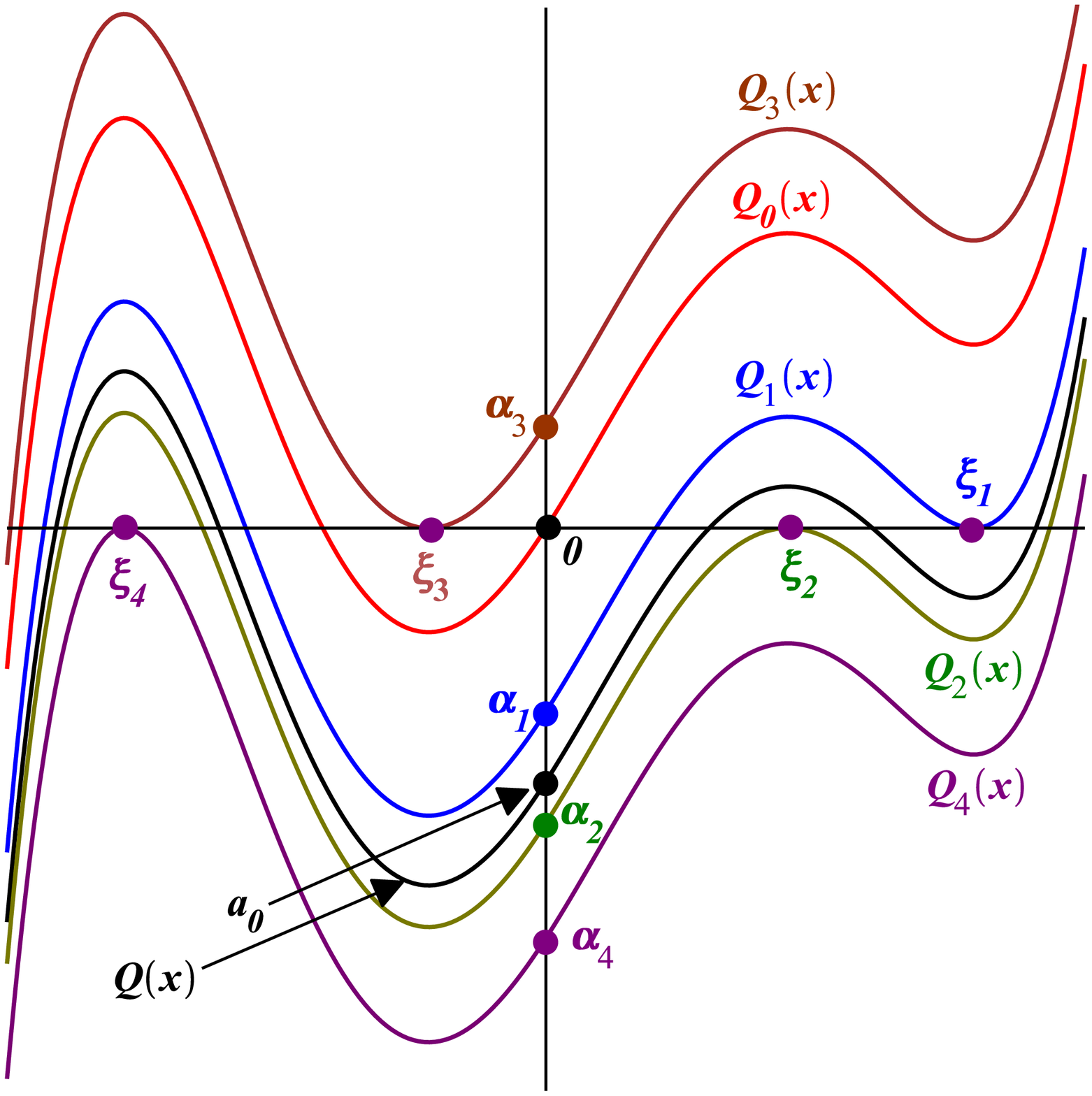} & \includegraphics[width=67mm]{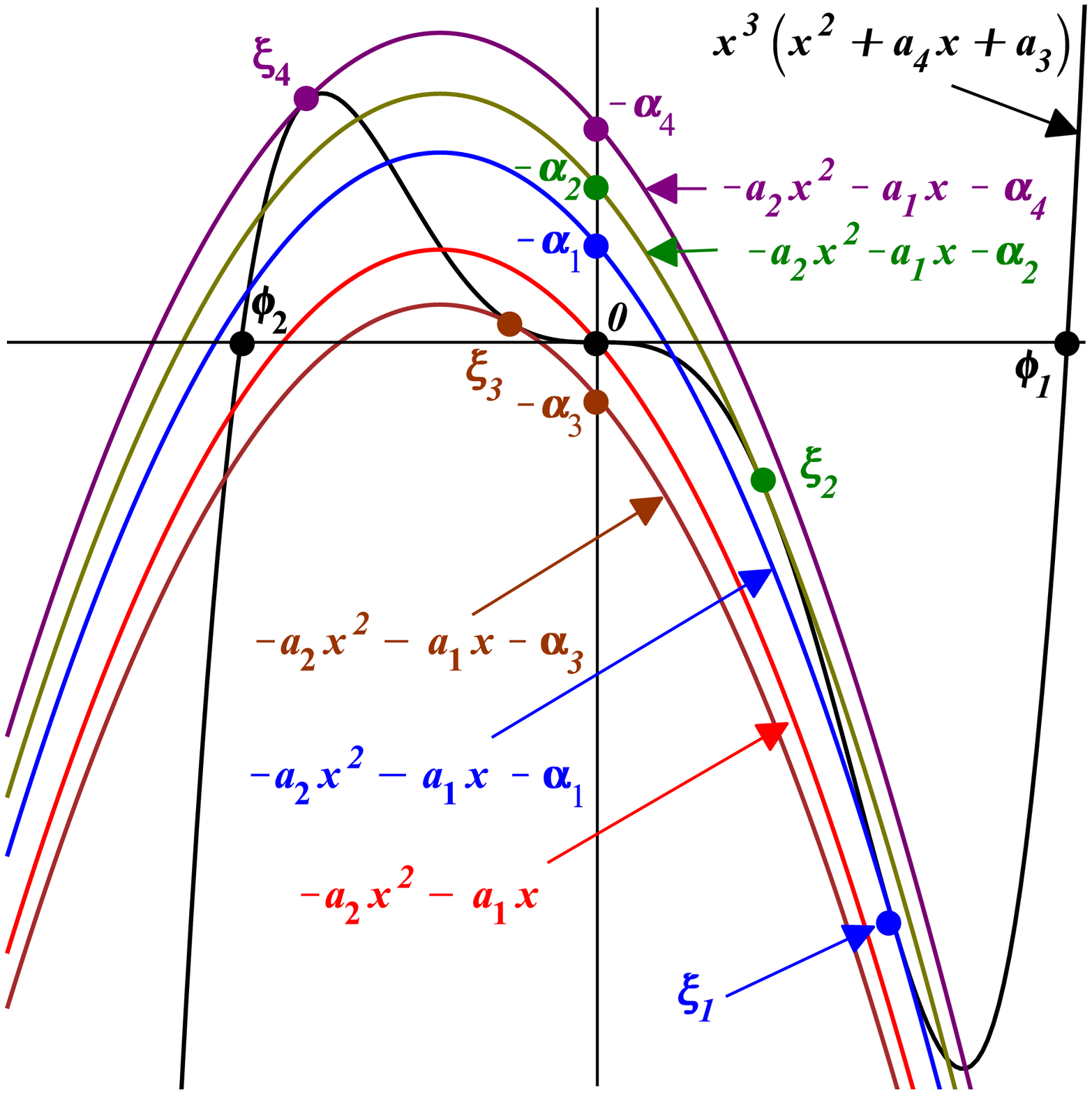} \\
{\scriptsize {\bf Figure 1a}} &  {\scriptsize {\bf Figure 1b}} \\
& \\
\multicolumn{1}{c}{\begin{minipage}{18em}
\scriptsize
\vskip-.15cm
A sample quintic $Q(x)$ with five real roots, the separatrix quintic $Q_0(x)$, and the four privileged quintics $Q_i(x)$ for which the abscissa is tangent to their graph at the stationary point $\xi_i$, i.e. the stationary points of $Q(x)$ are stationary points and also double roots for $Q_i(x)$.
\end{minipage}}
& \multicolumn{1}{c}{\begin{minipage}{18em}
\scriptsize
\vskip-.15cm
The two ``components" of the privileged quintics $Q_i(x)$, namely, the quintic $x^3 (x^2 + a_4 x + a_3)$ (common for all of them) and the quadratic $a_2 x^2 + a_1 x + \alpha_i$. At the stationary points $\xi_i$, the curves $x^3 (x^2 + a_4 x + a_3)$ and  $- a_2 x^2 - a_1 x - \alpha_i$ have the same value and their tangents coincide.
\end{minipage}}
\\
\\
\end{tabular}
\end{center}

\n
Then, it suffices to order the determined $\alpha_i$'s into a set with increasing values, with the given $a_0$ in its place within this set, and hence find the isolation interval of each root of the given quintic equation in a straightforward manner. \\
For example, suppose that on Figure 1b, the quadratic $q_2(x) = -a_2 x^2 - a_1 x - a_0$ (not shown on the graph) has free term $a_0$ satisfying $-\alpha_3 < 0 < -\alpha_1 < -\alpha_2 < -a_0 < -\alpha_4$. That is, suppose that the parabola $q_2(x) = -a_2 x^2 - a_1 x - a_0$ is between the uppermost two parabolas on Figure 1b. Then, the roots $x_i$ of the quintic equation satisfy the following. There is a negative root $x_3$ which is between $\xi_4$ and the bigger of the smaller root of the quadratic equation $a_2 x^2 + a_1 x + \alpha_2 = 0$ and the smaller root $\phi_2$ of the quadratic equation $x^2 + a_4 x + a_3 = 0$.  There is another negative root, $x_2$, which is greater than $\xi_4$. There is also a positive root $x_1$ between $\xi_2$ and $\phi_1$ --- the bigger root of the quadratic equation $x^2 + a_4 x + a_3 = 0$. The remaining two roots are complex. \\
As another example, suppose that on Figure 1b one has $-\alpha_3 < 0 < -\alpha_1  < -a_0 < -\alpha_2 < -\alpha_4$, i.e. the parabola $q_2(x) = -a_2 x^2 - a_1 x - a_0$ is between the second and the third parabolas from the top. Then the isolation intervals of the roots of the quintic are as follows. As before, there is a negative root $x_5$ between $\xi_4$ and the bigger of the smaller root of the quadratic equation $a_2 x^2 + a_1 x + \alpha_2 = 0$ and the smaller root $\phi_2$ of the quadratic equation $x^2 + a_4 x + a_3 = 0$. Also as before, there is a negative root $x_4$ greater than $\xi_4$. Next, there is a positive root $x_3$ between the bigger root of $a_2 x^2 + a_1 x + \alpha_1 = 0$ and $\xi_2$. Another positive root, $x_2$, lies between $\xi_2$ and $\xi_1$. Finally, there is a third positive root $x_1$ between $\xi_1$ and $\phi_1$ --- the bigger root of the quadratic equation $x^2 + a_4 x + a_3 = 0$. \\
As a final example in the vein of Figure 1b, if one has $- \alpha_4 < - a_0$, i.e. the parabola $q_2(x) = -a_2 x^2 - a_1 x - a_0$ is above the uppermost parabola, then the quintic equation has a single positive root greater than the smaller of the bigger root of $a_2 x^2 + a_1 x + a_0 = 0$ and the bigger root $\phi_1$ of $x^2 + a_4 x + a_3 = 0$ and smaller than the larger of the latter two.
\n
As the goal of this work is to propose a method for localization of the roots of the quintic equation by solving quadratic equations only and without recourse to cubic and quartic equations, one cannot have knowledge of the stationary points $\xi_i$. Respectively, the $\alpha_i$'s cannot be known either. Hence, some residual indeterminacy will remain --- similar, but significantly reduced than the one associated with the Descartes' rule of signs. It will still be possible to determine the isolation intervals of some of the roots of the quintic equation by solving quadratic equations with coefficients taken from those of the quintic equation. The indeterminacy can be lifted by using the complete root classification of the quintic (see the Discussion in Section 5). \\
Consider again the auxiliary quartic equation (\ref{quartic}) for the stationary points $\xi_i$ of the quintic equation (\ref{quintic}). The idea is to determine, also by solving quadratic equations only, the number of real roots of this equation and their location. For that, one needs to determine the stationary points $\mu_i$ of the auxiliary quartic, that is, the curvature change points of the given quintic. These are the roots of the {\it first auxiliary cubic equation}
\b
\label{cubic1}
x^3 + \frac{3a}{4} x^2 + \frac{b}{2} x + \frac{c}{4} \,\, \equiv \,\, x^3 + \frac{3 a_4}{5} x^2 + \frac{3 a_3}{10} x + \frac{a_2}{10} = 0
\e
The discriminant of this equation is
\b
\Delta^{(3)}_1 = - \frac{1728}{25} a_2^2 - \frac{10368}{125} a_4 \left( \frac{4}{15} a_4^2 - a_3 \right) a_2 + \frac{3456}{125} a_3^2 \left( \frac{3}{10} a_4^2 - a_3 \right).
\e
Set $\Delta^{(3)}_1 = 0$ and consider the obtained as a quadratic equation for the unknown $a_2$ (with $a_3$ and $a_4$ treated as parameters). This will be the {\it third resolvent quadratic equation}
\b
\label{quadratic3}
x^2 + \frac{6 a_4}{5} \left( \frac{4 a_4^2}{15} - a_3 \right) x - \frac{2 a_3^2}{5} \left( \frac{3 a_4^2}{4} - a_3 \right) = 0.
\e
The roots $c_{1,2}$ of the third resolvent quadratic equation are:
\b
\label{tse}
c_{1,2}(a_3, a_4) = c_0(a_3, a_4) \pm \frac{\sqrt{2}}{25} \sqrt{(2 a_4^2 - 5 a_3)^3},
\e
where $c_0(a_3, a_4) = (3/5) a_4 a_3 - (4/25) a_4^3$. \\
If $c_2 \le a_2 \le c_1$, then the third resolvent quadratic equation will have two real roots and, hence, the discriminant $\Delta^{(3)}_1$ will be non-negative and the first auxiliary cubic equation (\ref{cubic1}) will have three real roots. Otherwise, (\ref{cubic1}) will have only one real root.
Therefore, if $c_2 \le a_2 \le c_1$, then the auxiliary quartic equation (\ref{quartic}) will have 0, or 2, or 4 real roots $\xi_i$, that is, the number of stationary points of the quintic equation could be 0, or 2, or 4 and, thus, the number of its real roots would be 1, or 3, or 5. If however, $a_2 \notin [c_2, c_1]$, then the number of roots of the auxiliary quartic equation (\ref{quartic}) would be either 0 or 2, i.e. the stationary points of the quintic would be either 0 or 2 and, hence, the number of its real roots would be either 1 or 3. \\
The auxiliary quartic equation (\ref{quartic}) can be approached in a manner similar to the one used for the quintic equation (\ref{quintic}) --- see [4] for the classification of the roots of the monic symbolic quartic equation in terms of the coefficients of the quartic. Firstly, one treats the free term $d$ as a parameter that can be varied. Therefore, one has a foliation of the $xy$-plane with a congruence of quartics differing by their free term only. Next, from within this congruence, one identifies the {\it separatrix quartic} $x^2 (x^2 + a x + b) + c x$ and the privileged quartics $x^2 (x^2 + a x + b) + c x + \delta_i$ (whose number is either 1 or 3) for which the abscissa is tangent to the graph at the stationary point $\mu_i$ of the auxiliary quartic, that is, the curvature change points of the quintic. Then, one ``splits" the quartics as follows (see [4] for details):
\b
\label{split2}
x^2 (x^2 + a x + b) = - c x - d,
\e
At the curvature change points $\mu_i$ of the quintic, the straight line $- c x - \delta_i$ (with $i > 0$) is tangent to the ``sub-quartic" $x^2 (x^2 + a x + b)$. \\
With the ``split" (\ref{split2}), one obtains a congruence of straight lines $-cx - d$ which foliate the $xy$-plane. \\
If one solves the auxiliary cubic equation (\ref{cubic1}) and determines the points $\mu_i$, from the auxiliary quartic (\ref{quartic}) equation, one immediately finds that
\b
\label{delta}
\delta_i = - \mu_i^4 - a \mu_i^3 - b \mu_i^2 - c \mu_i.
\e
Then, the isolation intervals of each root of the auxiliary quartic equation can be easily determined --- see [4]. However, the exact values of the roots $\mu_i$ will not be sought as only quadratic equations will be solved. 
\vskip.5cm
\n
The analysis of the two ``components" of the quintic, that is $x^3 q_1(x)$ and $q_2(x)$, is absolutely straightforward. \\
The left-hand side $x^3 q_1(x)$ of (\ref{split}) has a root at $x = 0$ of order at least 3. Its other two roots are the roots $\phi_{1,2}$ of the quadratic equation $q_1(x) = x^2 + a_4 x + a_3 = 0$, namely:
\b
\phi_{1,2} = -\frac{1}{2} a^4 \pm \frac{1}{2} \sqrt{a_4^2 - 4 a_3}.
\e
These are real for $a_3 \le (1/4) a_4^2$. \\
The stationary points of $x^3 q_1(x)$ are: the saddle at $x = 0$ and the points
\b
\chi_{1,2} = -\frac{2}{5} a_4 \pm \frac{2}{5} \sqrt{a_4^2 - \frac{15}{4} a_3}.
\e
The latter two exist if $a_3 \le (4/15) a_4^2$. \\
At the stationary points $\chi_{1,2}$, the left-hand side $x^3 q_1(x)$ of (\ref{split}) takes the values:
\b
\chi_{1,2}^3 \,\, q_1(\chi_{1,2}) \equiv f_{1,2} \!\!\! & = & \!\!\! \frac{1}{3125} \left(- 2 a_4 \pm \sqrt{4 a_4^2 - 15 a_3}\right)^{\!\! 3} \!\!\left( -2 a_4^2 + 10 a_3 \pm a_4 \sqrt{4 a_4^2 -15 a_3} \right). \nonumber \\
&&
\e
The curvature change points of the left-hand side $x^3 q_1(x)$ of (\ref{split}) are:
\b
\sigma_{1,2} = - \frac{3}{10} a_4 \pm \frac{3}{10} \sqrt{a_4^2 - \frac{10}{3} a_3}.
\e
These are real for $a_3 \le (3/10) a_4^2$. \\
All possibilities for $x^3 q_1(x)$ of (\ref{split}) are shown on Figures 2a to 2o. 
\vskip.5cm
\n
Should $a_2 \ne 0$, the roots of the right-hand side $q_2(x)$ of (\ref{split}) are:
\b
\psi_{1,2} = -\frac{a_1}{2 a_2} \pm \frac{1}{2 a_2} \sqrt{a_1^2 - 4 a_0 a_2}.
\e
These are real for $a_1^2 - 4 a_0 a_2 \ge 0.$ \\
The right-hand side of (\ref{split}) has either one stationary point,
\b
\omega = - \frac{a_1}{2a_2},
\e
should $a_2 \ne 0$, or no stationary points [$q_2(x)$ is a straight line when $a_2 = 0$]. \\
At its stationary point $\omega$, the quadratic $q_2(x)$ takes the value
\b
q_2(\omega) \equiv g = \frac{a_1^2}{4 a_2} - a_0.
\e
All possibilities for $q_2(x)$ of (\ref{split}) are shown on Figures 3a to 3i. 
\vskip.5cm
\n
Altogether, there are $135$ possible cases: each of the 15 cases on Figures 2a to 2o for $x^3 q_1(x)$ with each of the nine cases for $q_2(x)$ on Figures 3a to 3i. But, as can be seen from the graphs, not all cases are qualitatively different and, in addition, the existing symmetry reduces further the number of qualitatively different possibilities. \\
In Section 4, one of these cases will be analyzed fully. The remaining ones could be done in a similar manner.

\vskip.5cm
\subparagraph{\hskip-.6cm 3 \hskip0.2cm Application of the Method \vskip.5cm}
\vskip.3cm
\n
\hskip-1cm First, for any given $a_3$ and $a_4$ [that is, whatever the ``sub-quintic" $x^5 + a_4 x^4 + a_3 x^3$ is], the coefficient $a_2$ of the quadratic $- a_2 x^2 - a_1 x - a_0$ ``selects" whether the full quintic (\ref{quintic}) will have one, or three, or five real roots (should $a_2 \in [c_2, c_1]$) or whether it will have one or three real roots only (should $a_2 \notin [c_2, c_1]$). As first step in the analysis, the value of $a_2$ relative to $c_1$ and $c_2$ should be determined. \\
Second, the coefficient $a_1$ ``determines" the exact number of roots of the auxiliary quartic equation (\ref{quartic}) (see [4] for details), that is, the exact number of stationary points of the quintic (\ref{quintic}). By solving cubic equations, the isolation interval of each real root of the auxiliary quartic equation can be determined [4]. If cubic equations are not to be solved and the analysis is done by solving quadratic equations only, then there would be some residual indeterminacy [4]. This steps necessitates the determination, following [4], of the number of roots of the auxiliary quartic equation (\ref{quartic}), that is, the number of stationary points of the quintic.   \\
Third, the free term $a_0$ ``selects" the position of the parabola $-a_2 x^2 - a_1 x - a_0$ among the congruence of parabolas and hence, the exact number of roots of the quintic (\ref{quintic}). When the stationary points of the quintic are known and, thus, the $\alpha_i$'s are known, then, by putting the elements of the set of all $\alpha_i$'s in increasing order and by determining the place of $a_0$ in that set, the isolation interval of each real root of the quintic can be immediately found. The end-points of the isolation intervals would be the stationary points of the quintic and the roots of the resolvent quadratic equations. Otherwise, when only quadratic equations are solved and, hence, the $\alpha_i$'s cannot be explicitly known, then, as already mentioned, there would be some residual indeterminacy. What could be determined in this case are either isolation intervals or root clustering intervals (the latter containing 1 or 3 real roots, or containing  0 or 2 real roots, or, very rarely, containing 0, or 2, or 4 real roots) with end-points given by the roots of the resolvent quadratic equations.

\begin{center}
\begin{tabular}{ccccc}
\includegraphics[width=27mm]{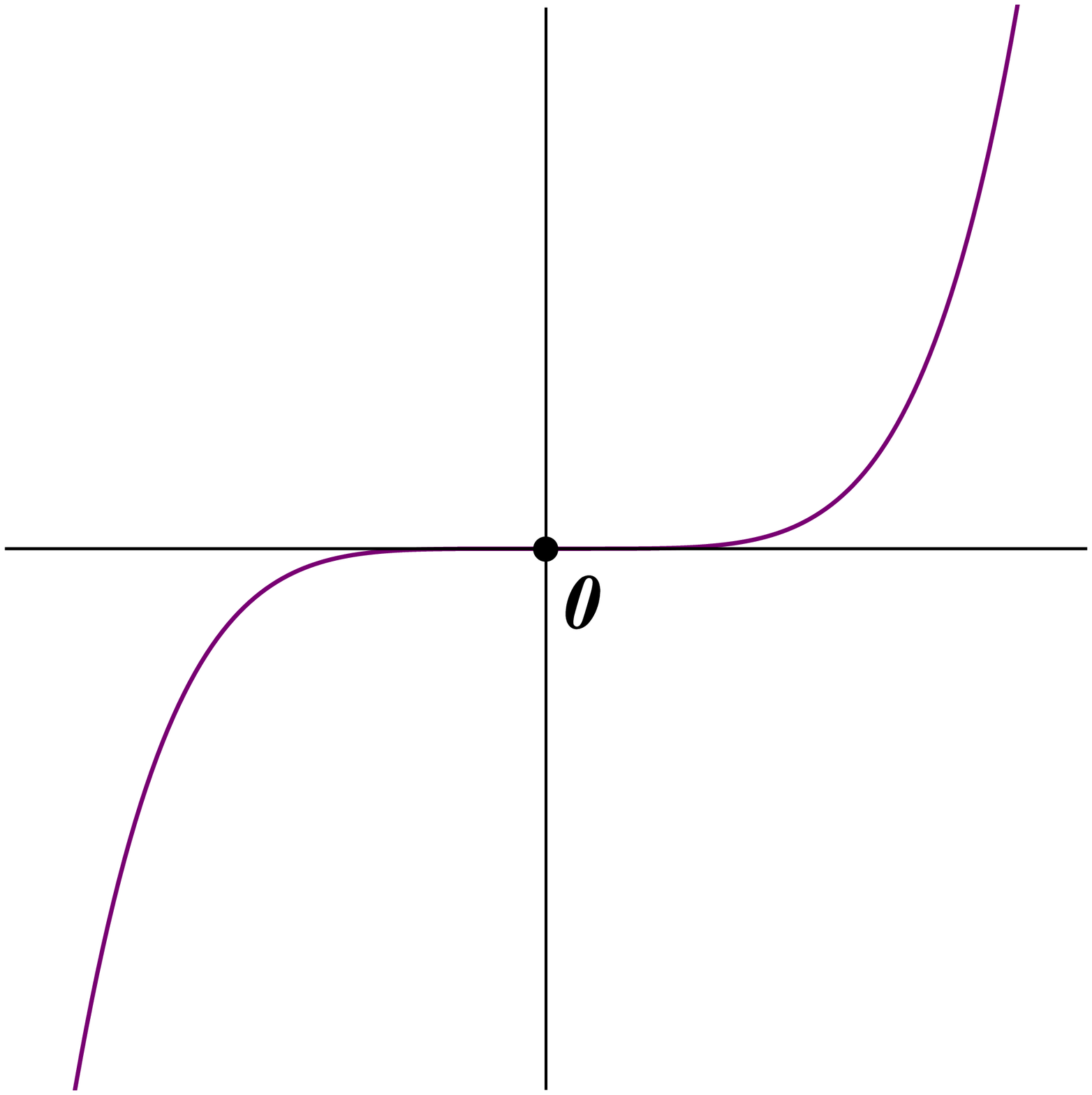} & \includegraphics[width=27mm]{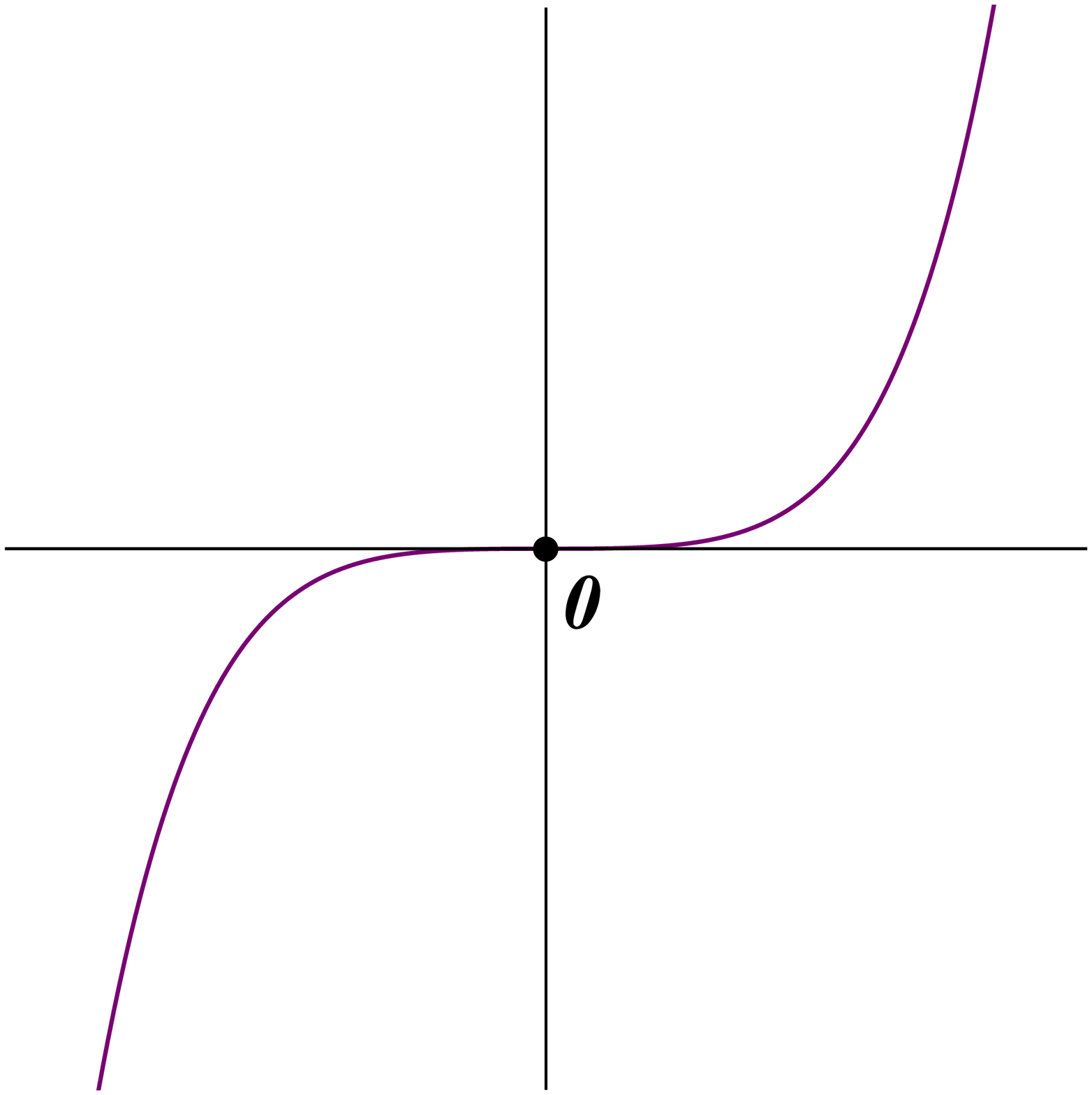} & \includegraphics[width=27mm]{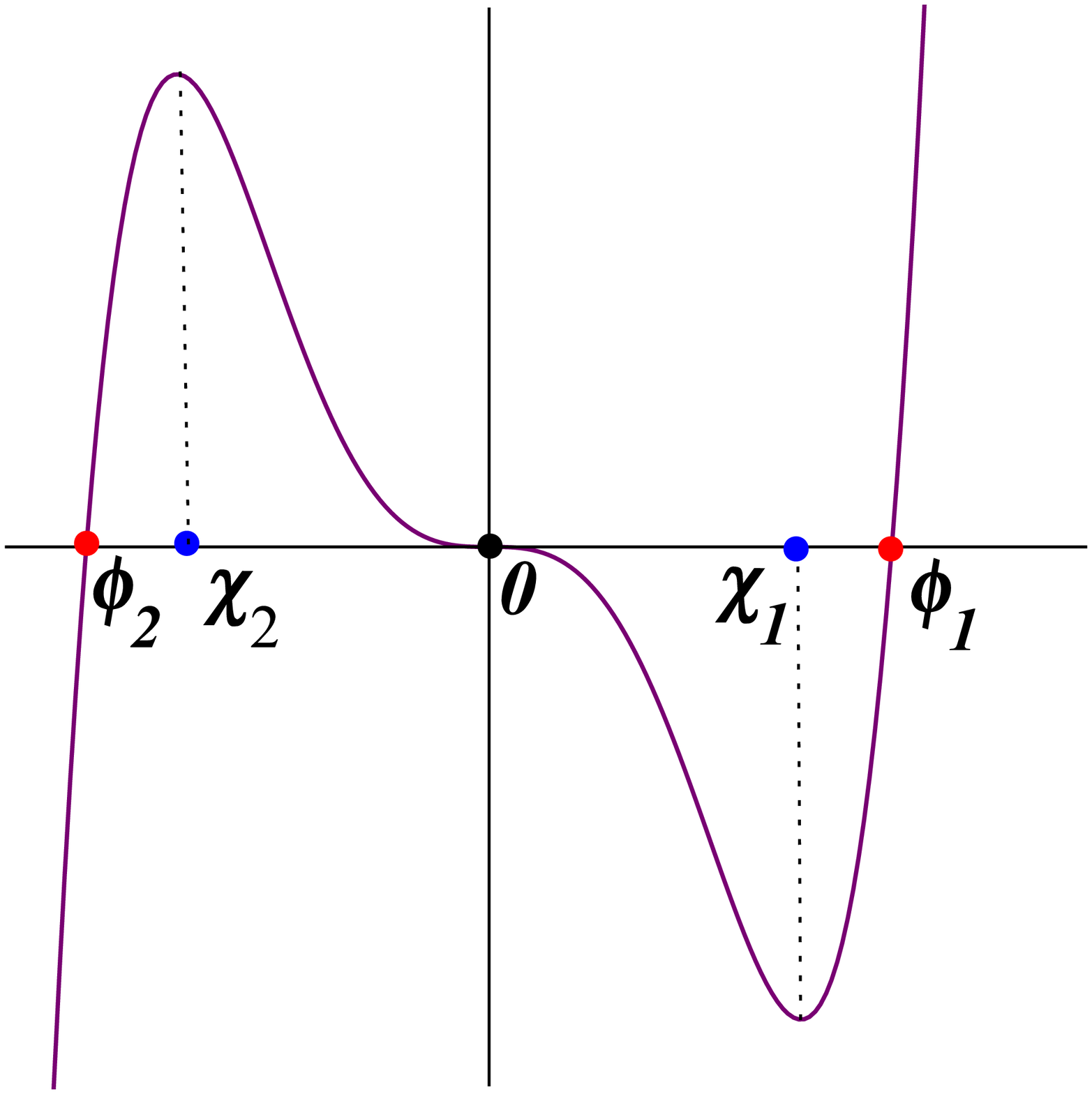} & \includegraphics[width=27mm]{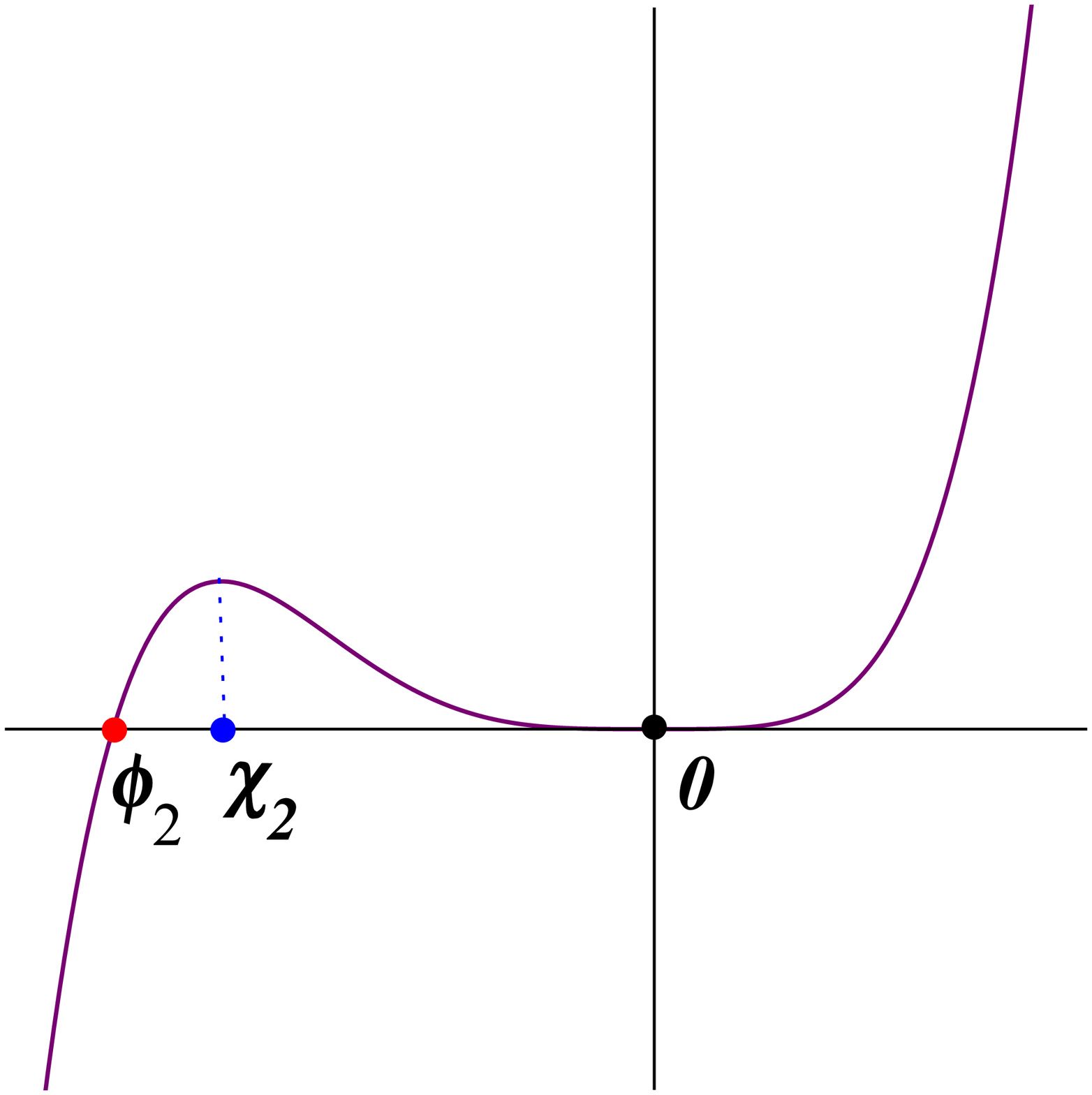} & \includegraphics[width=27mm]{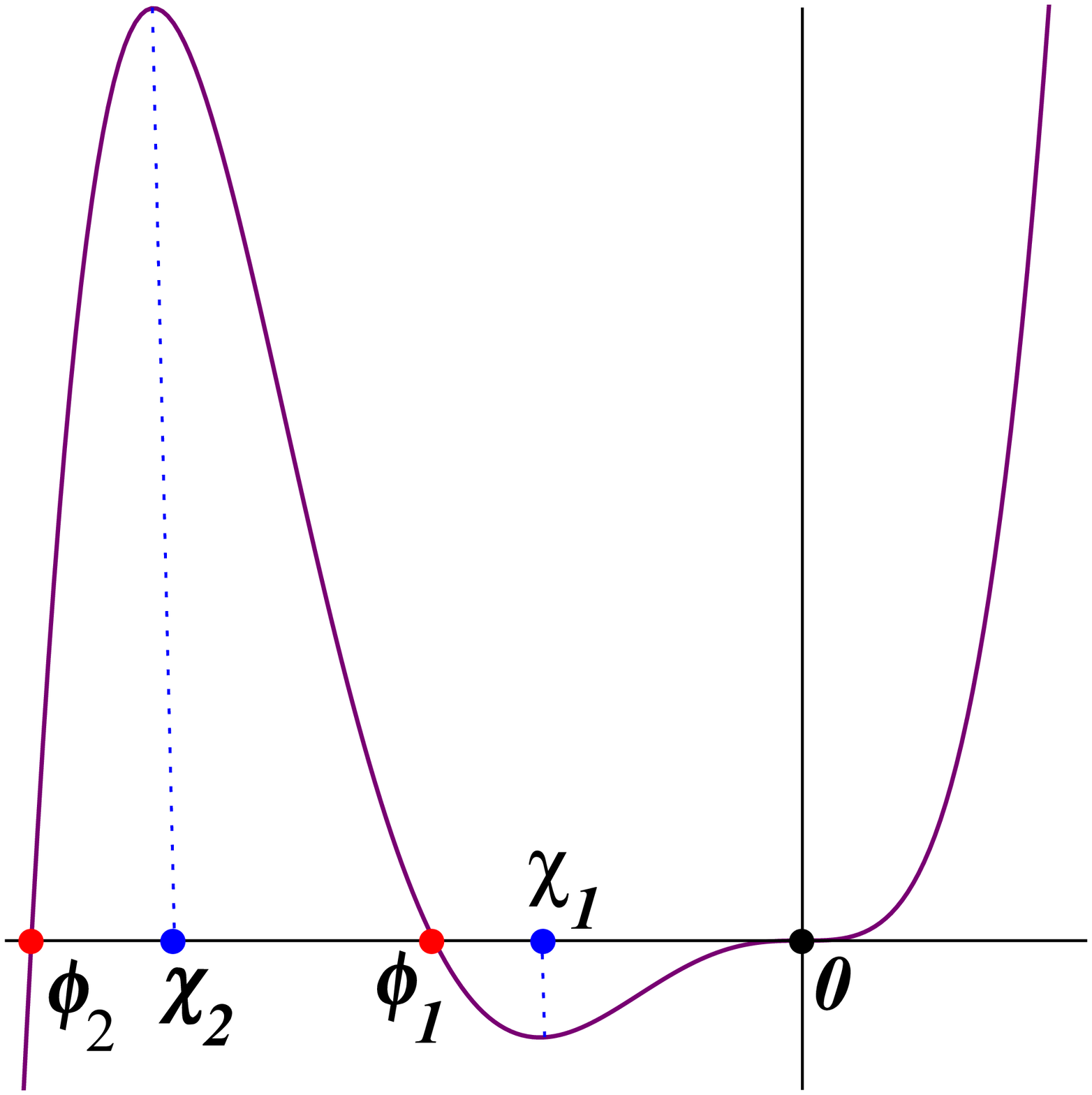}\\
{\scriptsize {\bf Figure 2a}} &  {\scriptsize {\bf Figure 2b}} & {\scriptsize {\bf Figure 2c}} &  {\scriptsize {\bf Figure 2d}} & {\scriptsize {\bf Figure 2e}} \\
& & & & \\
\multicolumn{1}{c}{\begin{minipage}{7em}
\scriptsize
\begin{center}
\vskip-4.35cm
$\bm{a_4 = 0}$ and $\bm{a_3 = 0}$
\end{center}
The point $x = 0$ is a quintuple root of $x^3 q_1(x) = 0$. The saddle at $x = 0$ is the only stationary point.
\end{minipage}}
& \multicolumn{1}{c}{\begin{minipage}{7em}
\scriptsize
\begin{center}
\vskip-3.7cm
$\bm{a_4 = 0}$ and $\bm{a_3 > 0}$
\end{center}
The point $x = 0$ is a triple root of $x^3 q_1(x) = 0$. The other two roots of $x^3 q_1(x) = 0$ are complex. The saddle at $x = 0$ is the only stationary point.
\end{minipage}}
& \multicolumn{1}{c}{\begin{minipage}{7em}
\scriptsize
\begin{center}
\vskip-1.75cm
$\bm{a_4 = 0}$ and $\bm{a_3 < 0}$
\end{center}
The point $x = 0$ is a triple root of $x^3 q_1(x) = 0$. The other two roots of $x^3 q_1(x) = 0$ are $\phi_{1,2} = \pm \sqrt{-a_3}$. The stationary points of $x^3 q_1(x)$ are: $\chi_{1,2} = \pm (1/5) \sqrt{- 15 a_3}$ (a local maximum at $\chi_2 < 0$ and a local minimum at $\chi_1 > 0$) and the saddle at $x = 0$.
\end{minipage}}
& \multicolumn{1}{c}{\begin{minipage}{7em}
\scriptsize
\begin{center}
\vskip-2.05cm
$\bm{a_4 > 0}$ and $\bm{a_3 = 0}$
\end{center}
\vskip-.01cm
The  point $x = 0$ is a quadruple root of $x^3 q_1(x) = 0$ (including $\phi_1 = 0$). The other root of $x^3 q_1(x) = 0$ is $\phi_2 = - a_4$. The stationary points are a local minimum at $x = 0$ (including $\chi_1 = 0$) and a local maximum at $\chi_2 = -(4/5) a_4$.
\end{minipage}}
& \multicolumn{1}{c}{\begin{minipage}{7em}
\scriptsize
\begin{center}
\vskip-.35cm
$\bm{a_4 > 0}$ and $\bm{0 < a_3 \le (1/4) a_4^2}$
\end{center}
\vskip-.35cm
The  point $x = 0$ is a triple root of $x^3 q_1(x) = 0$. The other two roots of $x^3 q_1(x) = 0$ are $\phi_{1,2} = -(1/2) a^4 \pm (1/2) \sqrt{a_4^2 - 4 a_3}$ --- both negative. The stationary points of $x^3 q_1(x)$ are: $\chi_{1,2} = -(2/5) a_4 \pm (1/5) \sqrt{4 a_4^2 - 15 a_3}$ (a local maximum at $\chi_2 < 0$ and a local minimum at $\chi_1 < 0$) and the saddle at $x = 0$.
\end{minipage}}
\\
\includegraphics[width=27mm]{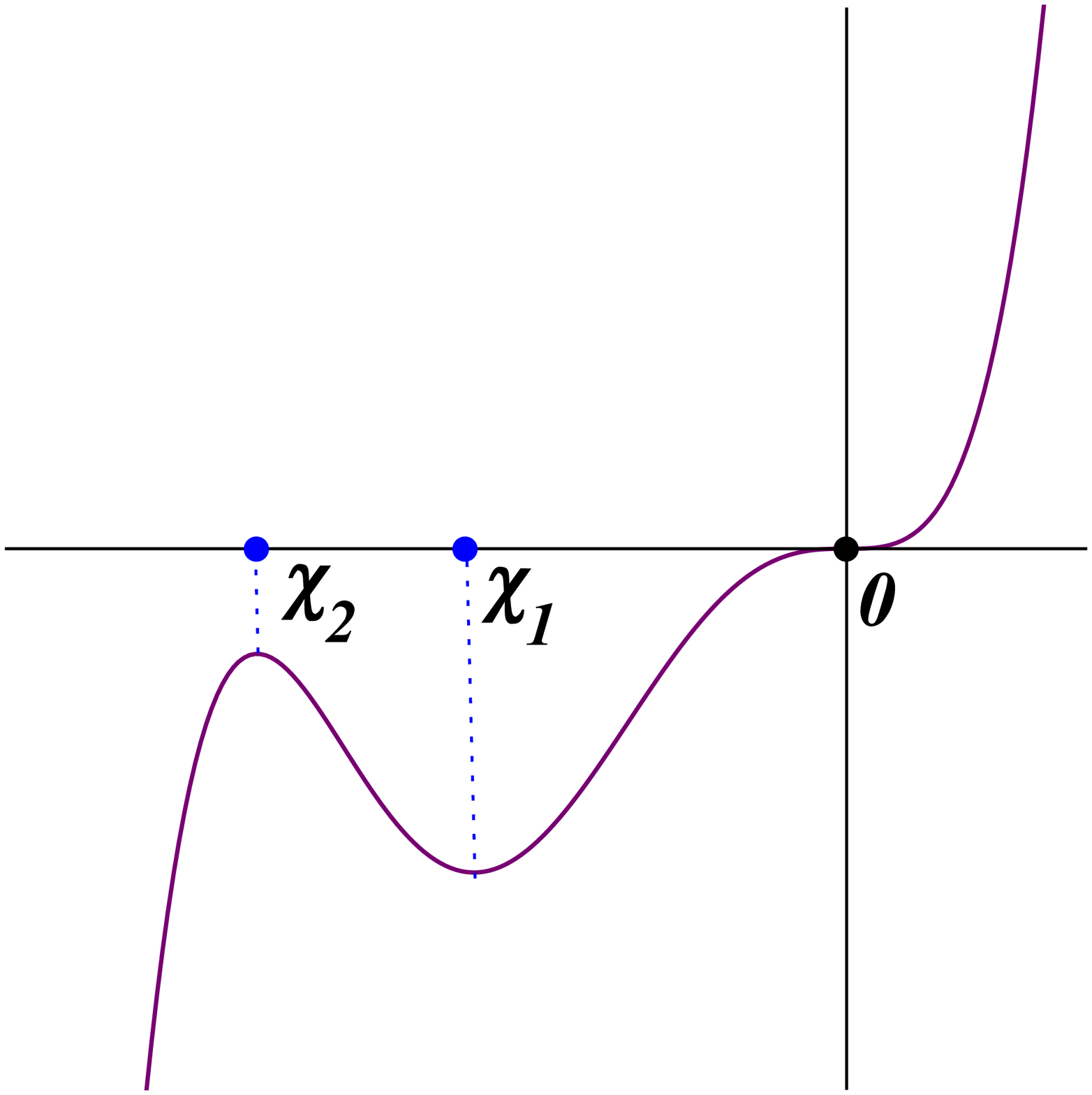} & \includegraphics[width=27mm]{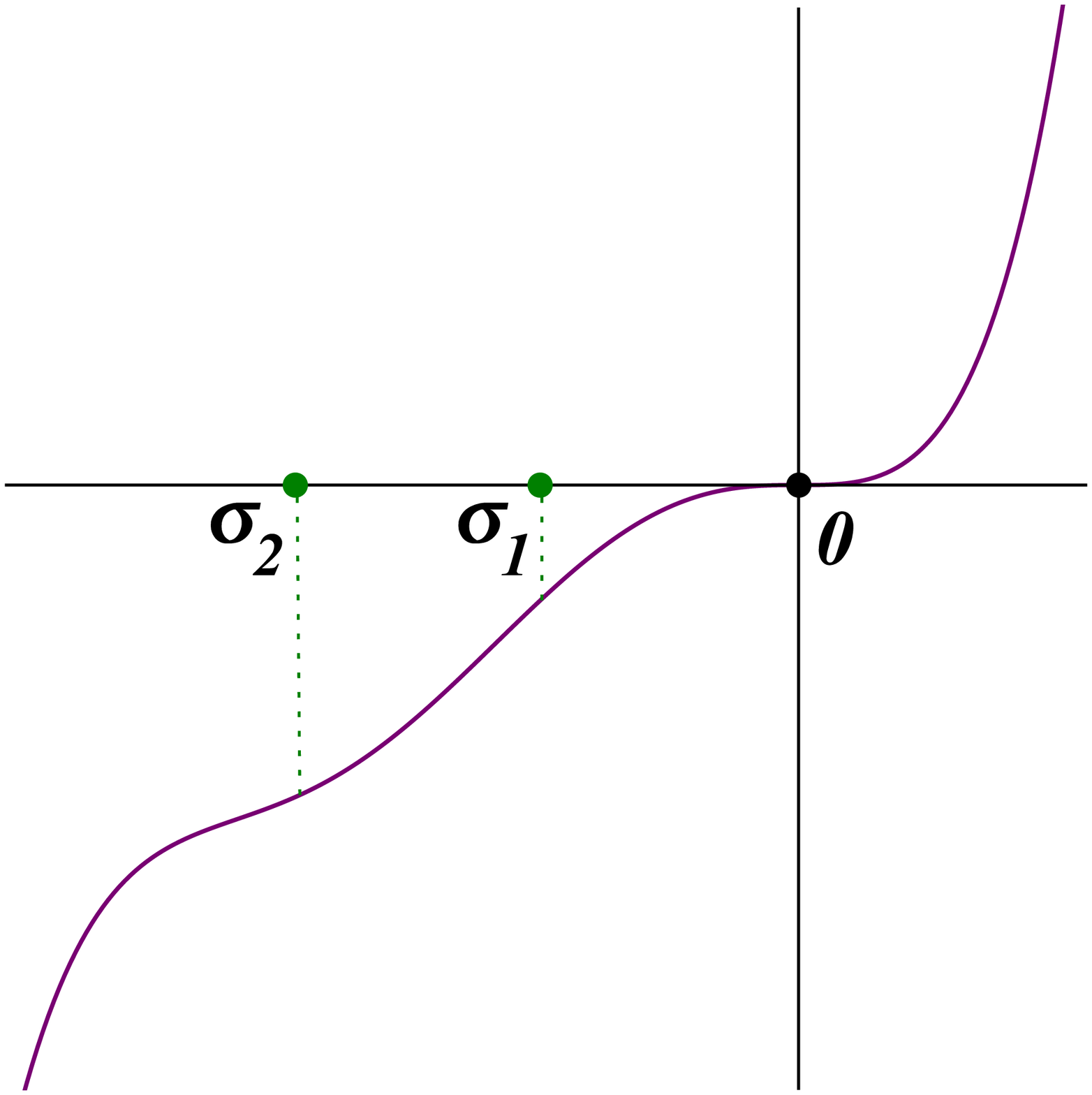} & \includegraphics[width=27mm]{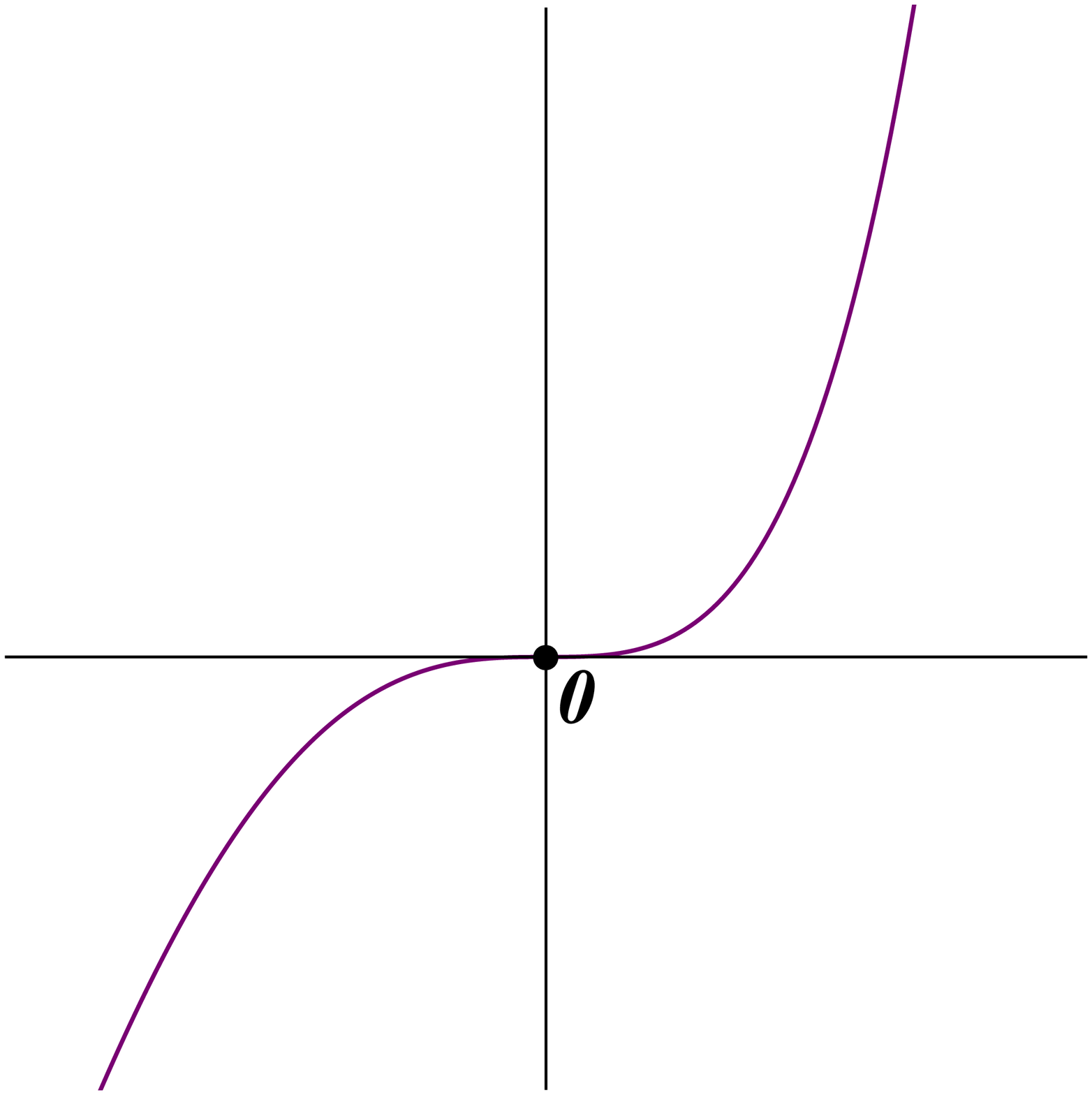} & \includegraphics[width=27mm]{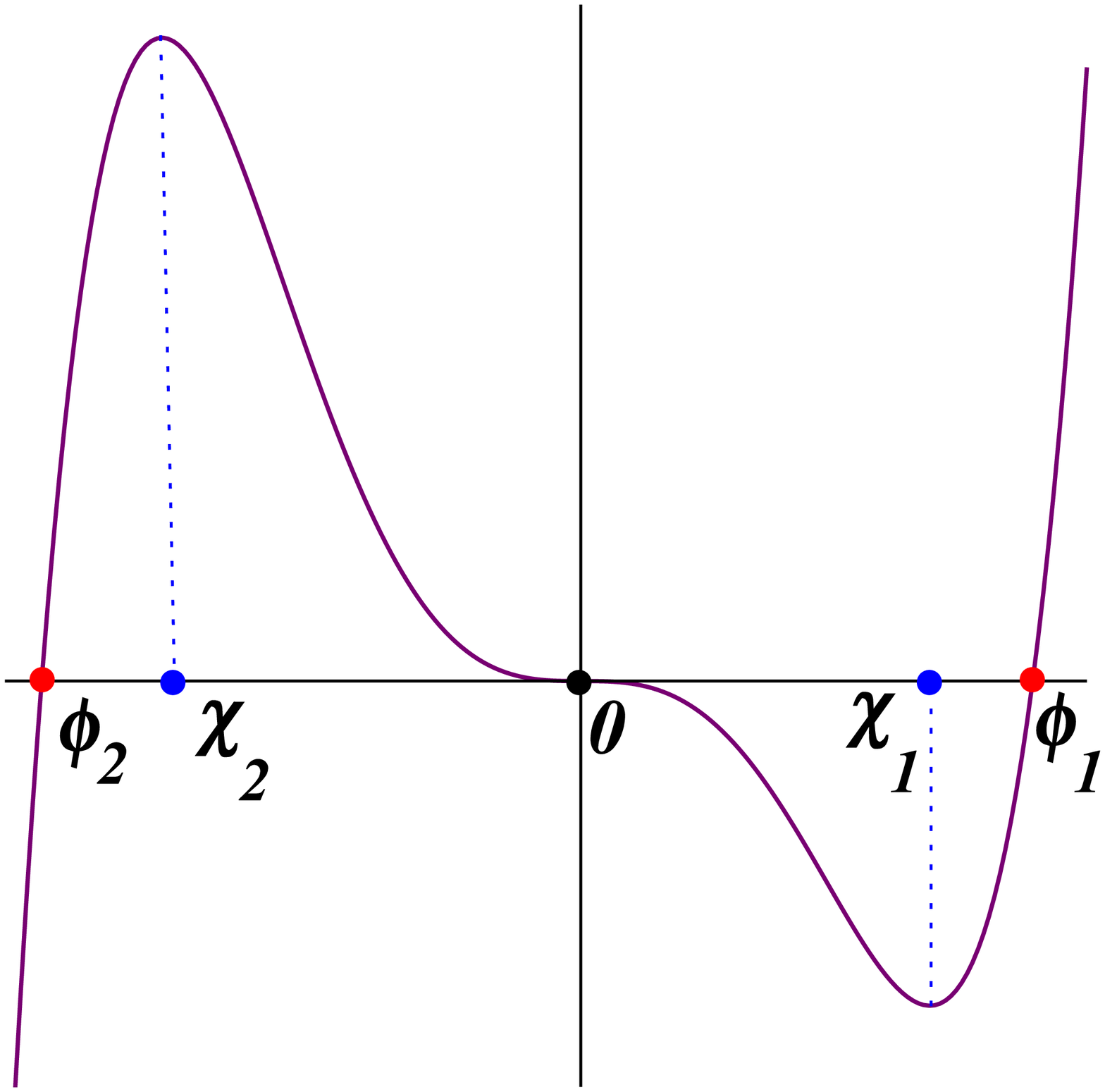} & \includegraphics[width=27mm]{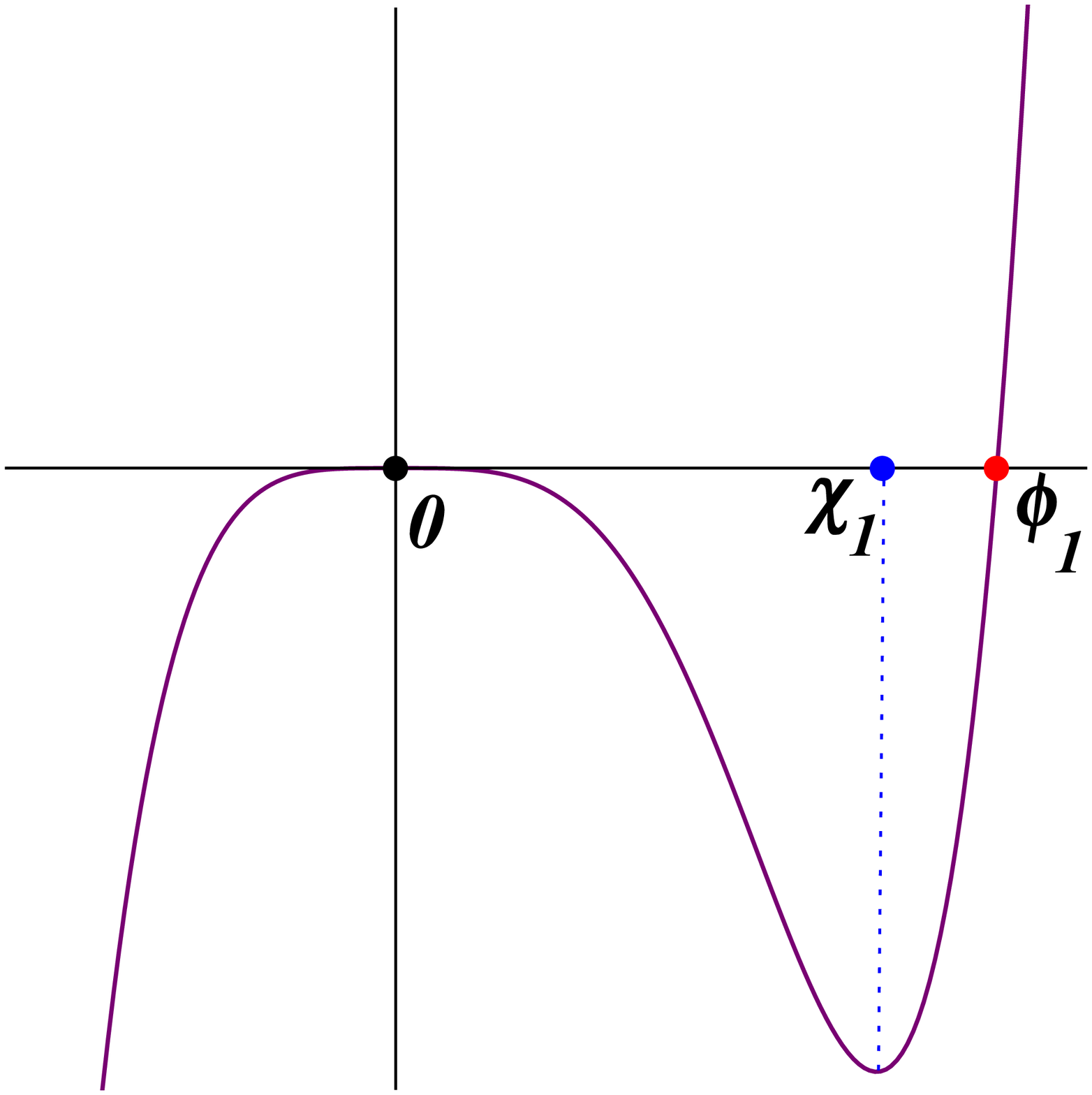}\\
{\scriptsize {\bf Figure 2f}} &  {\scriptsize {\bf Figure 2g}} & {\scriptsize {\bf Figure 2h}} &  {\scriptsize {\bf Figure 2i}} & {\scriptsize {\bf Figure 2j}} \\
& & & & \\
\multicolumn{1}{c}{\begin{minipage}{6.6em}
\scriptsize
\begin{center}
\vskip-.5cm
$\bm{a_4 > 0}$ and $\bm{(1/4) a_4^2 < a_3 \le (4/15) a_4^2}$
\end{center}
The  point $x = 0$ is a triple root of $x^3 q_1(x) = 0$. The other two roots of $x^3 q_1(x) = 0$ are complex. The quintic $x^3 q_1(x)$ still exhibits the two stationary points $\chi_{1,2} = -(2/5) a_4 \pm (1/5) \sqrt{4 a_4^2 - 15 a_3}$ (a local maximum at $\chi_2 < 0$ and a local minimum at $\chi_1 < 0$). In addition to them, $x^3 q_1(x)$ has a saddle at $x = 0$.
\end{minipage}}
& \multicolumn{1}{c}{\begin{minipage}{6.6em}
\scriptsize
\begin{center}
\vskip-.45cm
$\bm{a_4 > 0}$ and $\bm{(4/15) a_4^2 < a_3  \le (3/10) a_4^2}$
\end{center}
The  point $x = 0$ is a triple root of $x^3 q_1(x) = 0$. The other two roots of $x^3 q_1(x) = 0$ are complex. The quintic $x^3 q_1(x)$ has a saddle at $x = 0$ and no further stationary points ($\chi_{1,2}$ are both complex in this case). There are two curvature change points $\sigma_{1,2} = -(3/10) a_4 \pm (1/10) \sqrt{9 a_4^2 - 30 a_3}$ (both negative).
\end{minipage}}
& \multicolumn{1}{c}{\begin{minipage}{6.6em}
\scriptsize
\begin{center}
\vskip-2.70cm
$\bm{a_4 > 0}$ and $\bm{(3/10) a_4^2 < a_3}$
\end{center}
\vskip.35cm
The  point $x = 0$ is a triple root of $x^3 q_1(x) = 0$. The other two roots of $x^3 q_1(x) = 0$ are complex. The quintic $x^3 q_1(x)$ has a saddle at $x = 0$ and no further stationary points ($\chi_{1,2}$ are both complex in this case). 
\end{minipage}}
& \multicolumn{1}{c}{\begin{minipage}{6.6em}
\scriptsize
\begin{center}
\vskip-.3cm
$\bm{a_4 > 0}$ and $\bm{a_3 < 0}$
\end{center}
The  point $x = 0$ is a triple root of $x^3 q_1(x) = 0$. The other two roots of $x^3 q_1(x) = 0$ are $\phi_{1,2} = -(1/2) a^4 \pm (1/2) \sqrt{a_4^2 - 4 a_3}$ ($\phi_1$ is positive and $\phi_2$ --- negative, with $|\phi_2| > |\phi_1|$). The stationary points of $x^3 q_1(x)$ are: $\chi_{1,2} = -(2/5) a_4 \pm (1/5) \sqrt{4 a_4^2 - 15 a_3}$ (a local maximum at $\chi_2 < 0$ and a local minimum at $\chi_1 > 0$) and the saddle at $x = 0$.
\end{minipage}}
& \multicolumn{1}{c}{\begin{minipage}{6.6em}
\scriptsize
\begin{center}
\vskip-2.3cm
$\bm{a_4 < 0}$ and $\bm{a_3 = 0}$
\end{center}
\vskip0.65cm
The  point $x = 0$ is a quadruple root of $x^3 q_1(x) = 0$ (including $\phi_2 = 0$). The other root of $x^3 q_1(x) = 0$ is $\phi_1 = - a_4$.
The stationary points are a local maximum at $x = 0$ (including $\chi_2 = 0$) and a local minimum at $\chi_1 = -(4/5) a_4$.
\end{minipage}}
\\
\end{tabular}
\end{center}

\begin{center}
\begin{tabular}{ccccc}
\includegraphics[width=27mm]{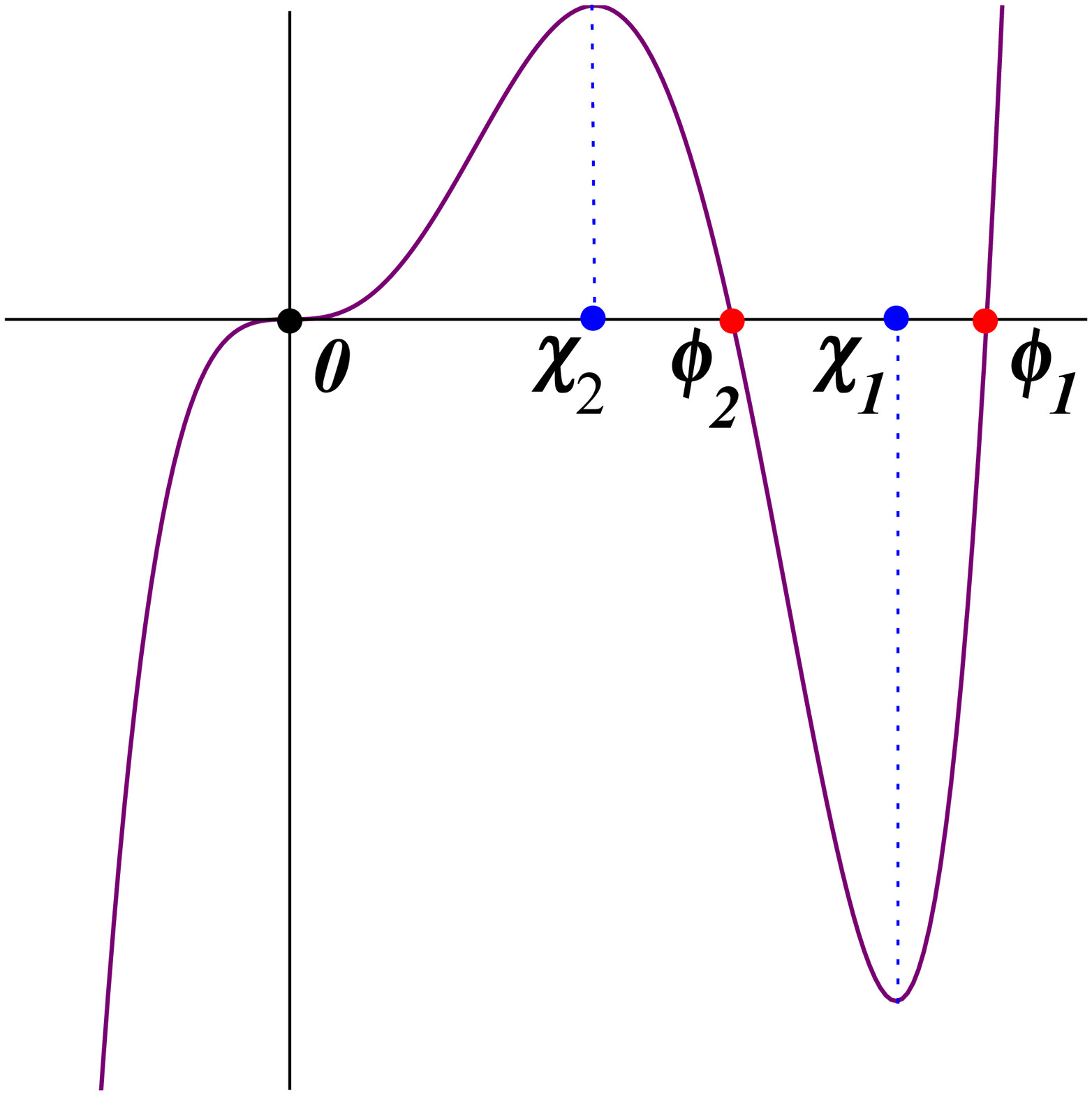} & \includegraphics[width=27mm]{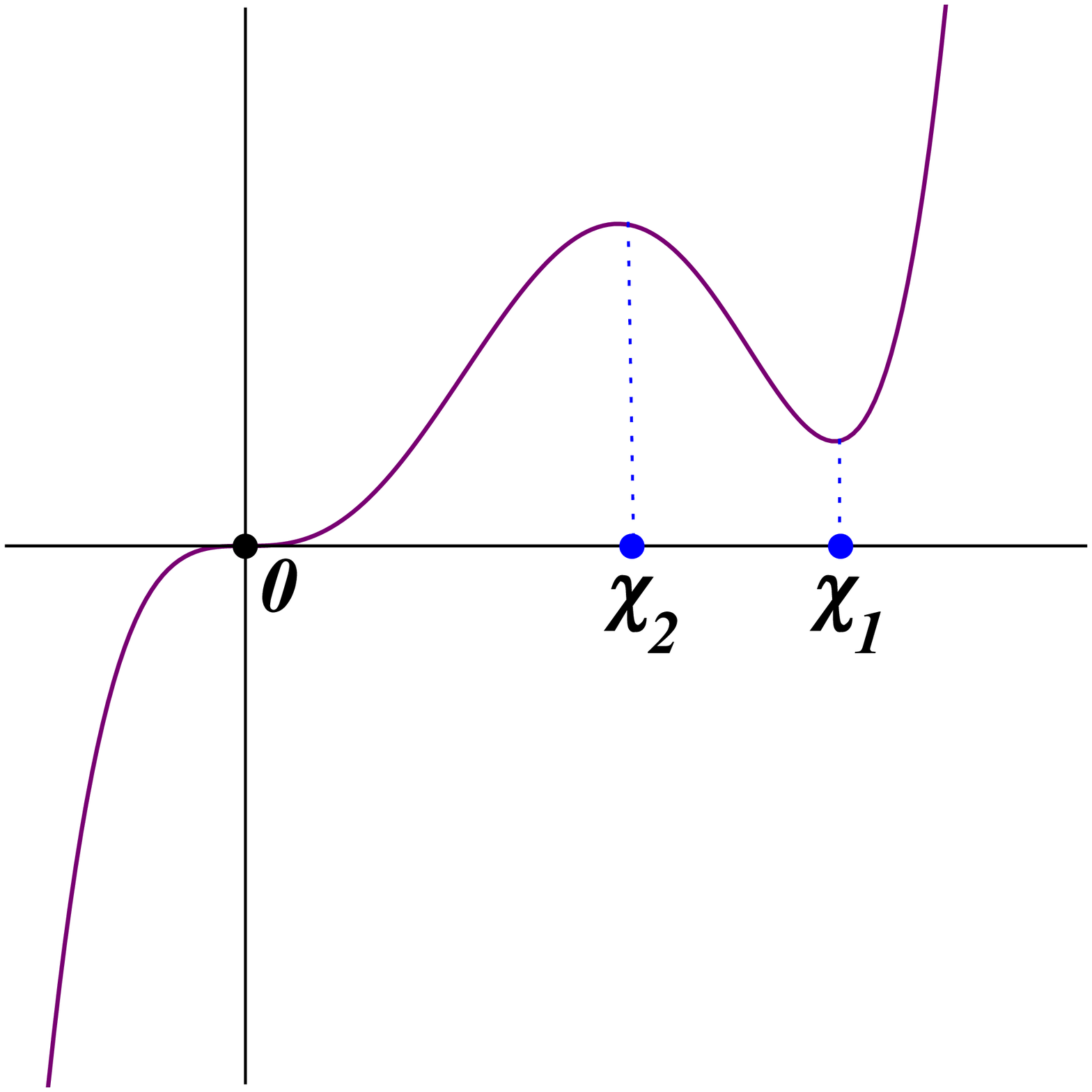} & \includegraphics[width=27mm]{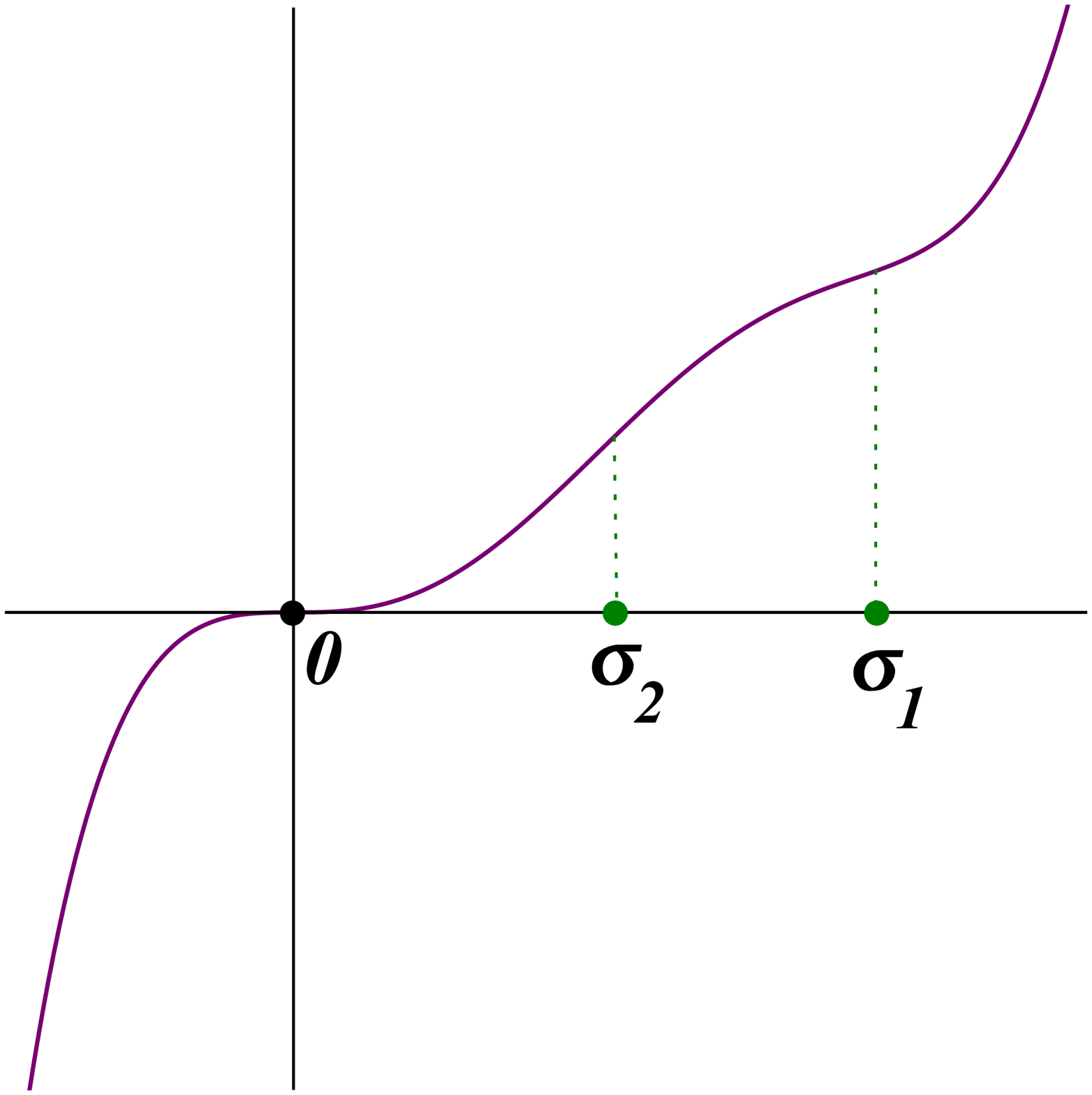} & \includegraphics[width=27mm]{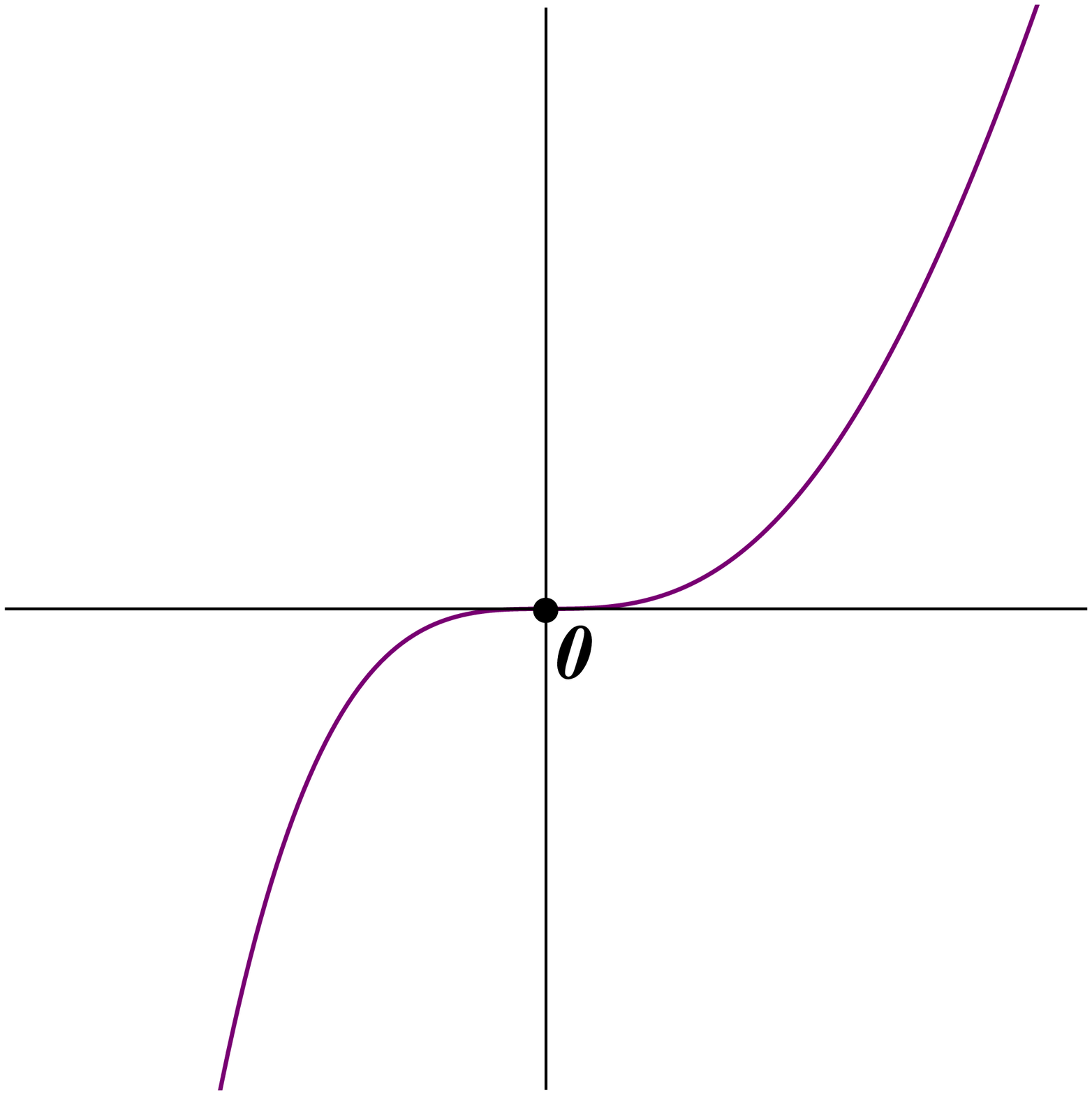} & \includegraphics[width=27mm]{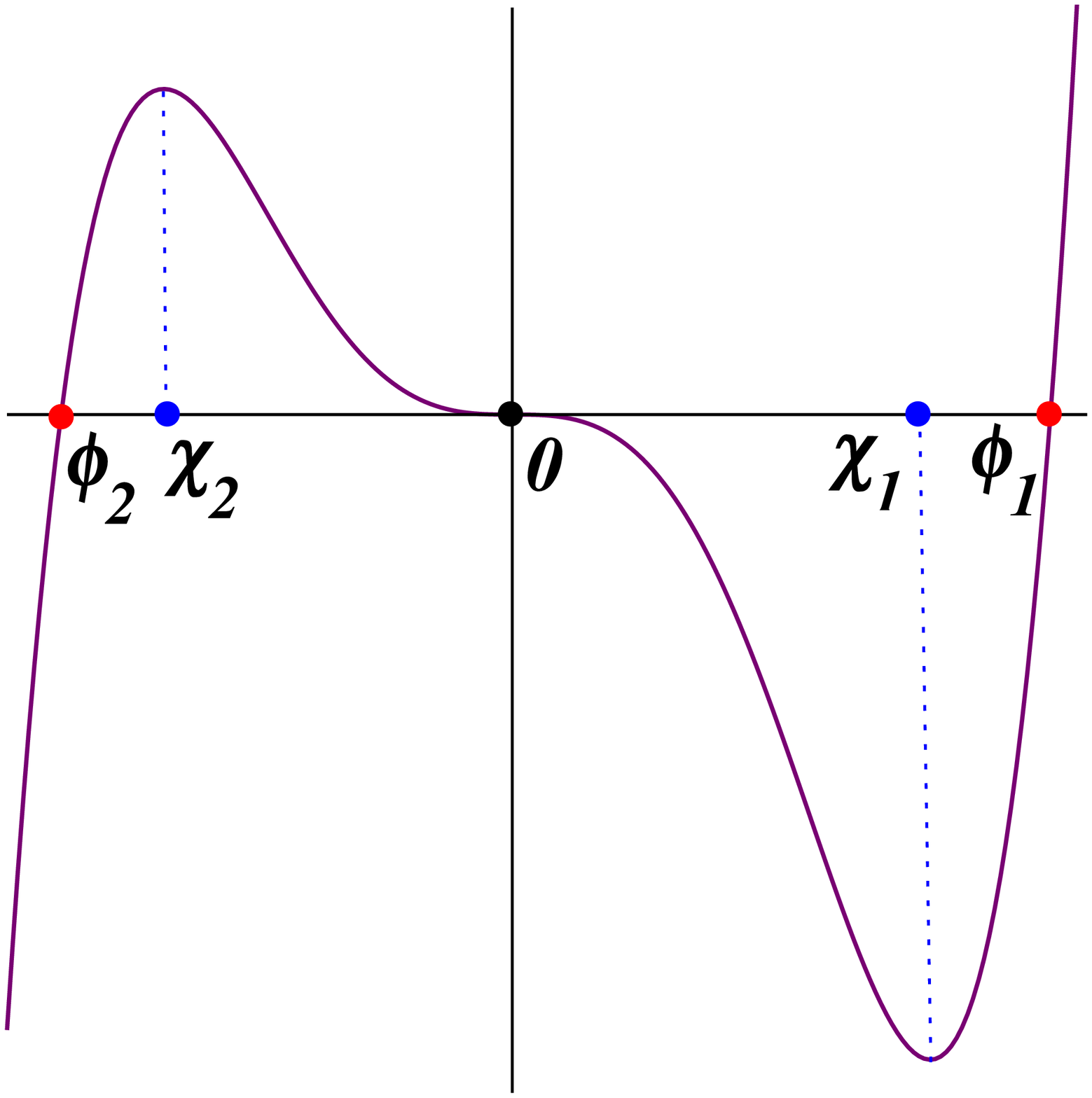}\\
{\scriptsize {\bf Figure 2k}} &  {\scriptsize {\bf Figure 2l}} & {\scriptsize {\bf Figure 2m}} &  {\scriptsize {\bf Figure 2n}} & {\scriptsize {\bf Figure 2o}} \\
& & & & \\
\multicolumn{1}{c}{\begin{minipage}{6.6em}
\scriptsize
\begin{center}
\vskip-.3cm
$\bm{a_4 < 0}$ and $\bm{0 < a_3 \le (1/4) a_4^2}$
\end{center}
\vskip.1cm
The  point $x = 0$ is a triple root of $x^3 q_1(x) = 0$. The other two roots of $x^3 q_1(x) = 0$ are $\phi_{1,2} = -(1/2) a^4 \pm (1/2) \sqrt{a_4^2 - 4 a_3}$ --- both positive. The quintic $x^3 q_1(x)$ still exhibits the two stationary points $\chi_{1,2} = -(2/5) a_4 \pm (1/5) \sqrt{4 a_4^2 - 15 a_3}$ (a local maximum at $\chi_2 > 0$ and a local minimum at $\chi_1 > 0$). In addition to them, $x^3 q_1(x)$ has a saddle at $x = 0$.
\end{minipage}}
& \multicolumn{1}{c}{\begin{minipage}{6.6em}
\scriptsize
\begin{center}
\vskip-1.2cm
$\bm{a_4 < 0}$ and $\bm{(1/4) a_4^2 < a_3 \le (4/15) a_4^2}$
\end{center}
\vskip-.2cm
The  point $x = 0$ is a triple root of $x^3 q_1(x) = 0$. The other two roots of $x^3 q_1(x) = 0$ are complex. The quintic $x^3 q_1(x)$ still exhibits the two stationary points $\chi_{1,2} = -(2/5) a_4 \pm (1/5) \sqrt{4 a_4^2 - 15 a_3}$ (a local maximum at $\chi_2 > 0$ and a local minimum at $\chi_1 > 0$). In addition to them, $x^3 q_1(x)$ has a saddle at $x = 0$.
\end{minipage}}
& \multicolumn{1}{c}{\begin{minipage}{6.6em}
\scriptsize
\begin{center}
\vskip-1.1cm
$\bm{a_4 < 0}$ and $\bm{(4/15) a_4^2 < a_3  \le (3/10) a_4^2}$
\end{center}
\vskip-.2cm
The  point $x = 0$ is a triple root of $x^3 q_1(x) = 0$. The other two roots of $x^3 q_1(x) = 0$ are complex. The quintic $x^3 q_1(x)$ has a saddle at $x = 0$ and no further stationary points ($\chi_{1,2}$ are both complex in this case). There are two curvature change points $\sigma_{1,2} = -(3/10) a_4 \pm (1/10) \sqrt{9 a_4^2 - 30 a_3}$ (both positive).
\end{minipage}}
& \multicolumn{1}{c}{\begin{minipage}{6.6em}
\scriptsize
\begin{center}
\vskip-3.3cm
$\bm{a_4 < 0}$ and $\bm{(3/10) a_4^2 < a_3}$
\end{center}
\vskip.1cm
The  point $x = 0$ is a triple root of $x^3 q_1(x) = 0$. The other two roots $\phi_{1,2}$ of $x^3 q_1(x) = 0$ are complex. The quintic $x^3 q_1(x)$ has a saddle at $x = 0$ and no further stationary points ($\chi_{1,2}$ are both complex in this case).
\end{minipage}}
& \multicolumn{1}{c}{\begin{minipage}{6.6em}
\scriptsize
\begin{center}
\vskip-.1cm
$\bm{a_4 < 0}$ and $\bm{a_3 < 0}$ 
\end{center}
\vskip.6cm
\vskip-.35cm
The  point $x = 0$ is a triple root of $x^3 q_1(x) = 0$. The other two roots of $x^3 q_1(x) = 0$ are $\phi_{1,2} = -(1/2) a^4 \pm (1/2) \sqrt{a_4^2 - 4 a_3}$ ($\phi_1$ is positive and $\phi_2$ --- negative, with $|\phi_2| < |\phi_1|$). The stationary points of $x^3 q_1(x)$ are: $\chi_{1,2} = -(2/5) a_4 \pm (1/5) \sqrt{4 a_4^2 - 15 a_3}$ (a local maximum at $\chi_2 < 0$ and a local minimum at $\chi_1 > 0$) and the saddle at $x = 0$.
\end{minipage}}
\\
\end{tabular}
\end{center}

\vskip.5cm
\subparagraph{\hskip-.6cm 4 \hskip0.2cm Example with $\bm{a_4 > 0, \,\, a_3 < 0, \,\, a_2 > 0, \,\, a_1 < 0}$ (Case of Figure 2i with Figure 3c) --- Full Analysis \vskip.5cm}

\hskip-1cm Consider first the quintic $Q^{(5)}_1(x) = x^5 + x^4 - 2 x^3 + (5/6) x^2 - (1/8) x + a_0$ on Figures 4a and 4b. As the coefficient of the quadratic term, $a_2 = 5/6$, is between the roots $c_2 = -3.71$ and $c_1 = 0.99$ of the third resolvent quadratic equation (\ref{tse}), the quintic can have zero, or two, or four stationary points --- these are the roots of the auxiliary quartic equation (\ref{quartic}) and, therefore, the number of real roots of the quintic could be one, or three, or five. The auxiliary quartic equation (\ref{quartic}) for this quintic is $x^4 + (4/5) x^3 - (6/5) x^2 + (1/3) x - 1/40 = 0$. The analysis on Figure 1.13 in [4] shows that the auxiliary quartic equation has four real roots (one of which is negative and the rest --- positive). Hence, depending on the value of the free term of the quintic, and also on the value of the linear term of the quintic, the number of real roots can be one, or three, or five. 

\begin{center}
\begin{tabular}{ccccc}
\includegraphics[width=27mm]{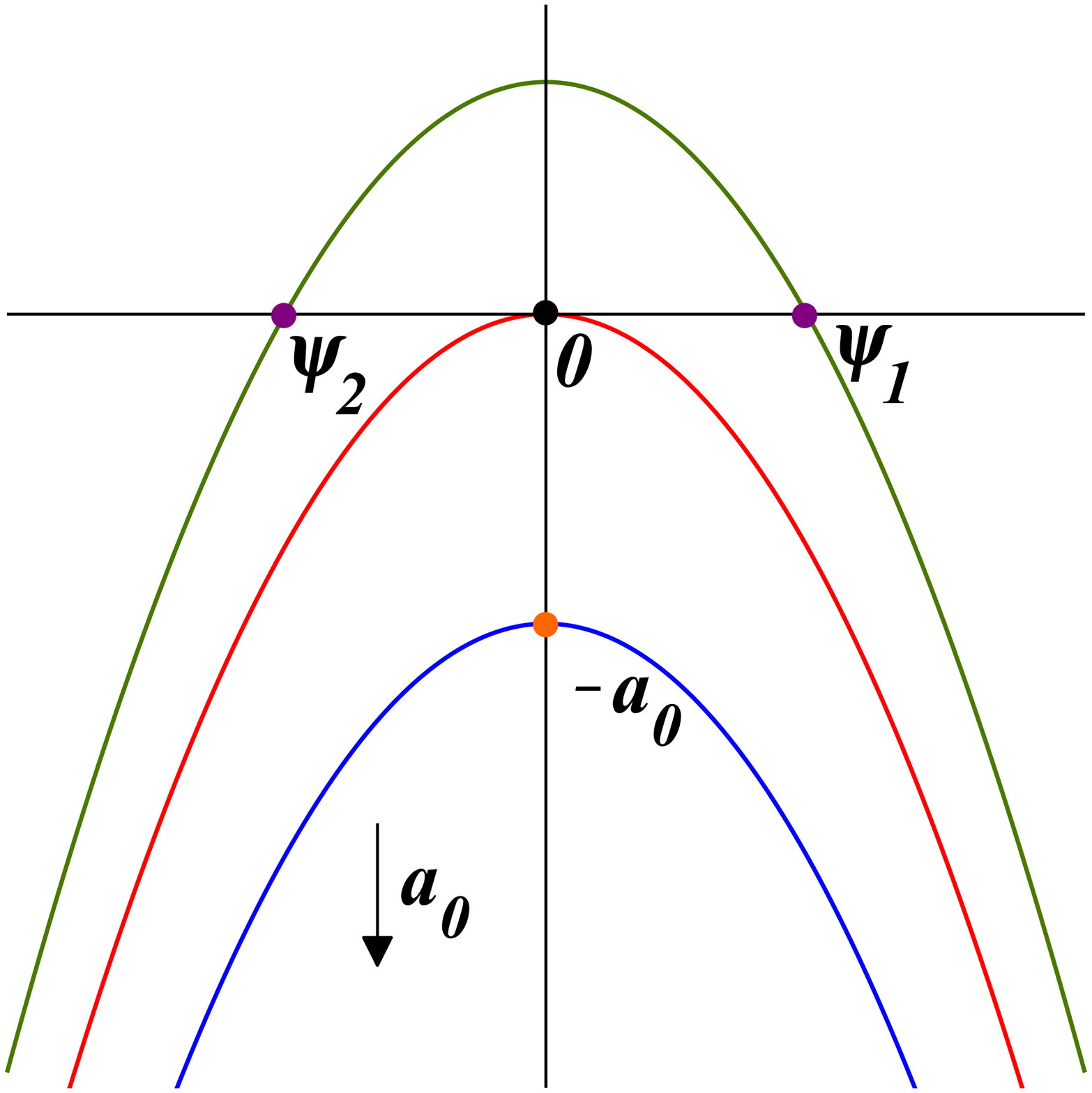} & \includegraphics[width=27mm]{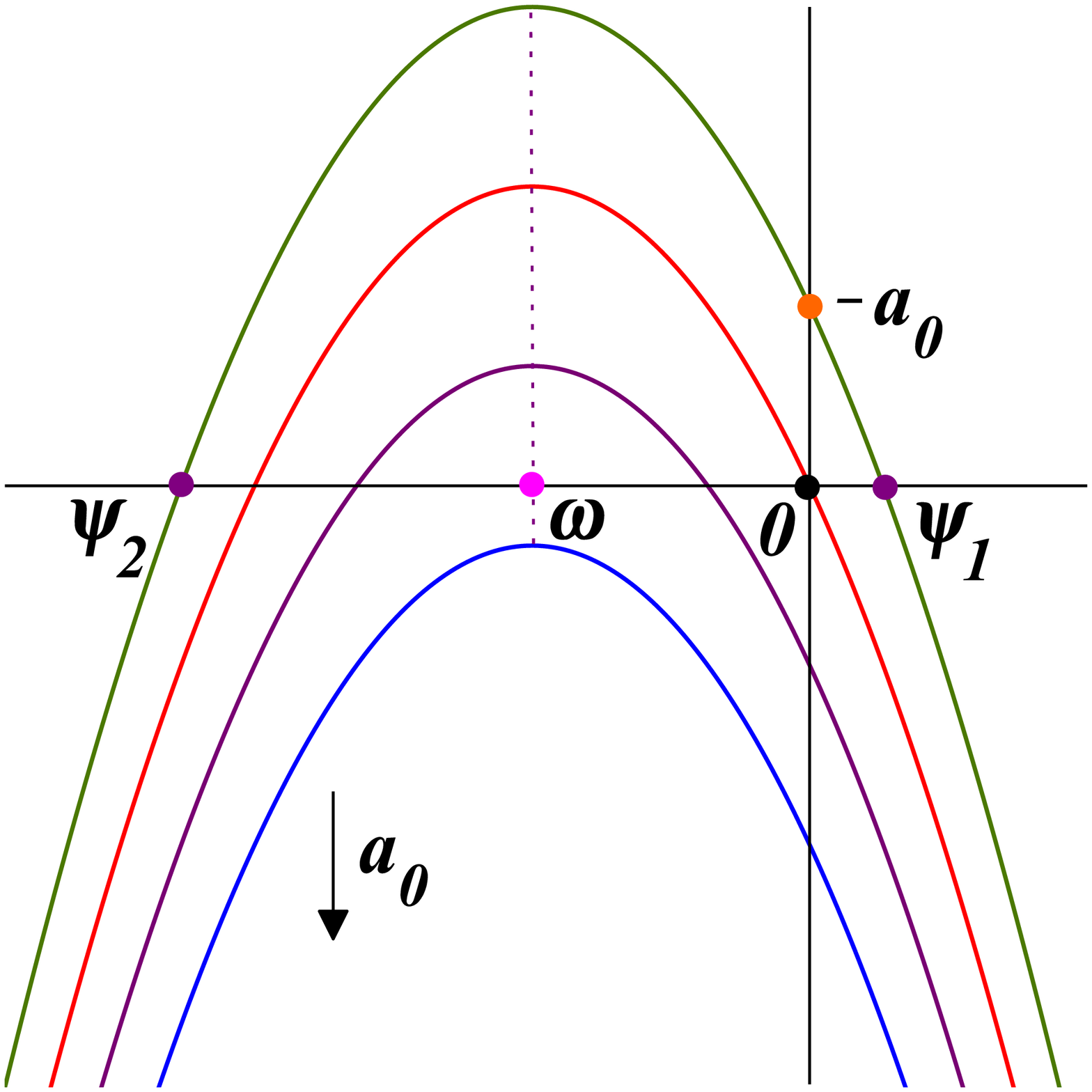} & \includegraphics[width=27mm]{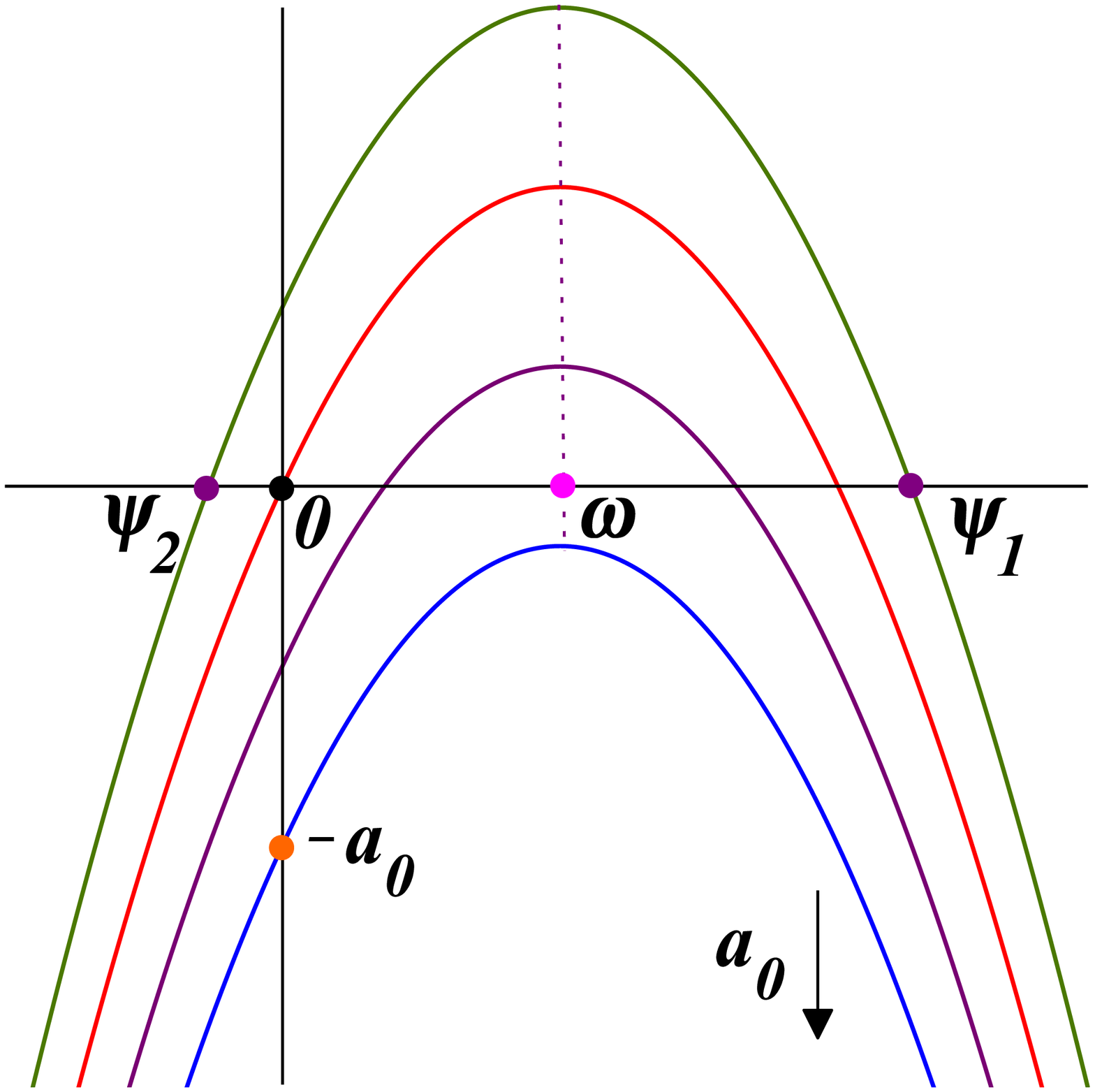} & \includegraphics[width=27mm]{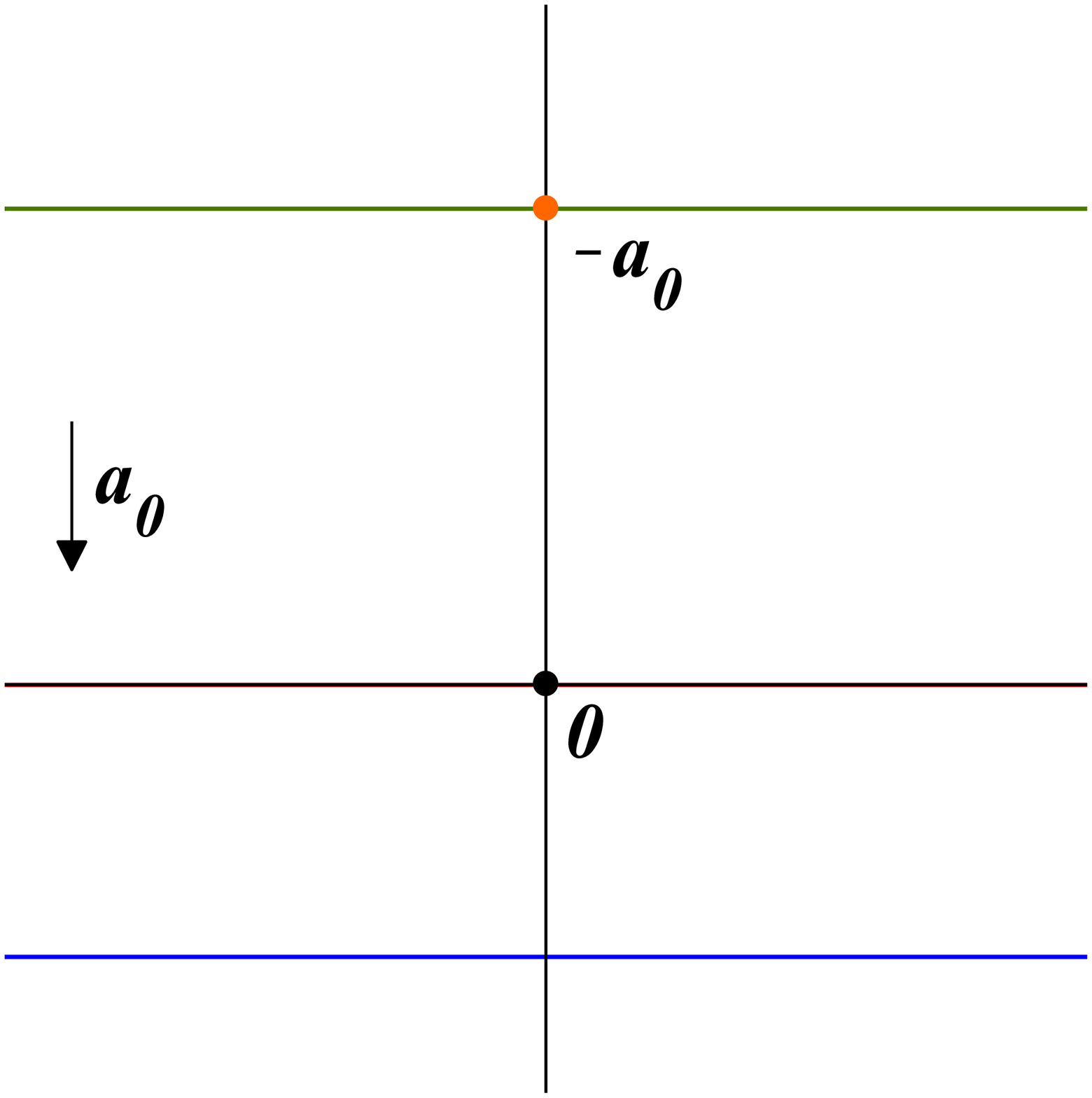} & \includegraphics[width=27mm]{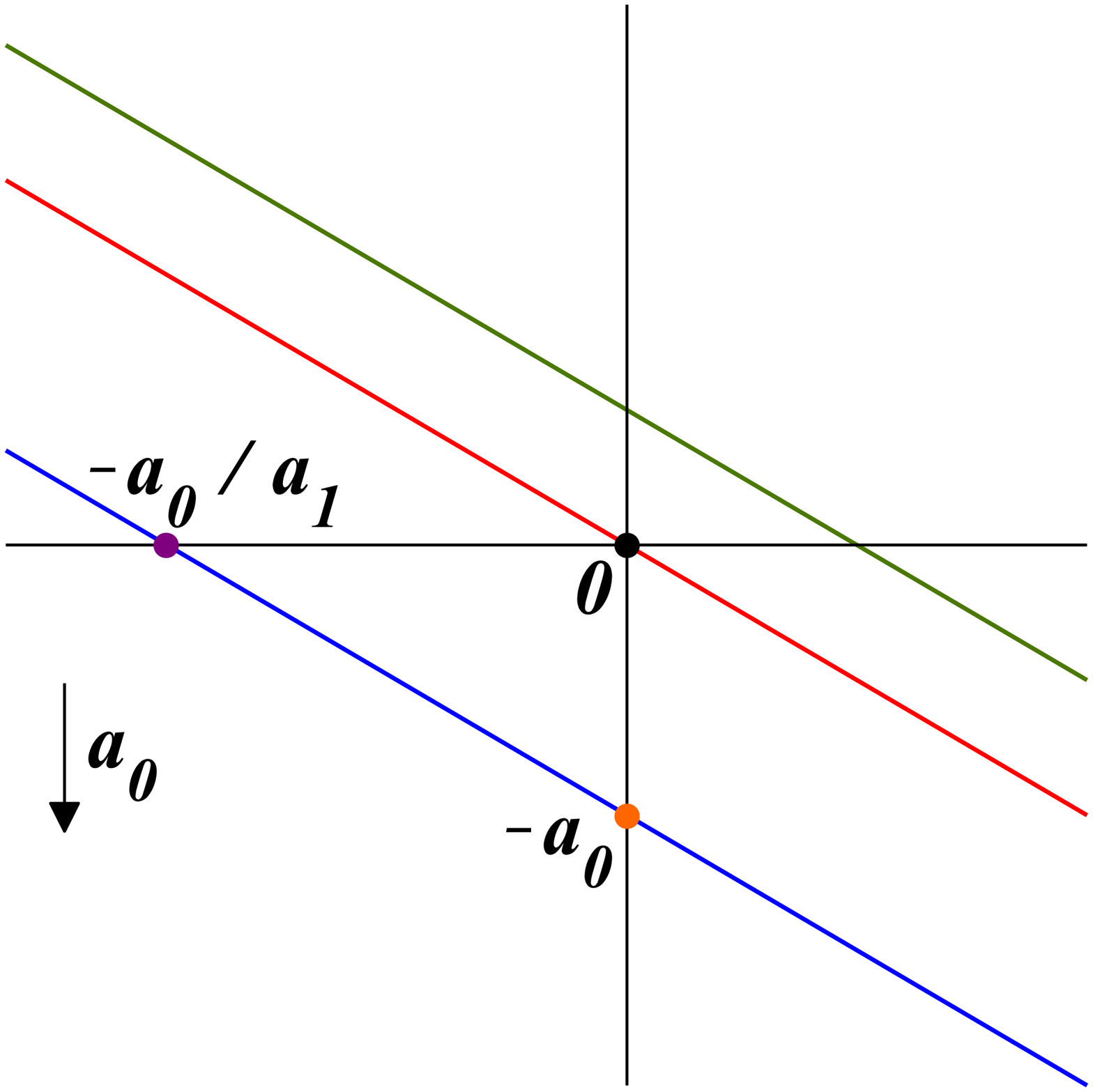}\\
{\scriptsize {\bf Figure 3a}} &  {\scriptsize {\bf Figure 3b}} & {\scriptsize {\bf Figure 3c}} &  {\scriptsize {\bf Figure 3d}} & {\scriptsize {\bf Figure 3e}} \\
& & & & \\
\multicolumn{1}{c}{\begin{minipage}{6.6em}
\scriptsize
\begin{center}
\vskip-2.5cm
$\bm{a_2 > 0}$ and $\bm{a_1 = 0}$ 
\end{center}
\vskip-.1cm
If $\bm{a_0 \le 0}$, the quadratic $q_2(x)$ has roots $\psi_{1,2} = \pm \sqrt{-a_0/a_2}$. If $\bm{a_0 > 0}$, the quadratic $q_2(x)$ has no real roots. In either case, there is a local maximum at $\omega = 0$.
\end{minipage}}
& \multicolumn{1}{c}{\begin{minipage}{6.6em}
\scriptsize
\begin{center}
\vskip-0.25cm
$\bm{a_2 > 0}$ and $\bm{a_1 > 0}$
\end{center}
\vskip-.1cm
If $\bm{a_0 \le a_1^2/(4 a_2)}$, the quadratic $q_2(x)$ has roots $\psi_{1,2} = - a_1/(2 a_2) \pm 1/(2 a_2)$ $\sqrt{a_1^2 - 4 a_0 a_2}$ ($\psi_2 < 0$, while $\psi_1$ could be negative, zero, or positive; $|\psi_2| > |\psi_1|$). If $\bm{a_0 > a_1^2/(4 a_2)}$, the quadratic $q_2(x)$ has no real roots. In either case, there is a local maximum at $\omega = -a_1/(2a_2) < 0$.
\end{minipage}}
& \multicolumn{1}{c}{\begin{minipage}{6.6em}
\scriptsize
\begin{center}
\vskip-0.25cm
$\bm{a_2 > 0}$ and $\bm{a_1 < 0}$
\end{center}
\vskip-.1cm
If $\bm{a_0 \le a_1^2/(4 a_2)}$, the quadratic $q_2(x)$ has roots $\psi_{1,2} = - a_1/(2 a_2) \pm 1/(2 a_2)$ $\sqrt{a_1^2 - 4 a_0 a_2}$ ($\psi_1 > 0$, while $\psi_2$ could be negative, zero, or positive; $|\psi_2| < \psi_1$). If $\bm{a_0 > a_1^2/(4 a_2)}$, the quadratic $q_2(x)$ has no real roots. In either case, there is a local maximum at $\omega = -a_1/(2a_2) > 0$.
\end{minipage}}
& \multicolumn{1}{c}{\begin{minipage}{6.6em}
\scriptsize
\begin{center}
\vskip-3.8cm
$\bm{a_2 = 0}$ and $\bm{a_1 = 0}$
\end{center}
\vskip-.15cm
The right-hand side $q_2(x)$ of (\ref{split}) is a straight line parallel to the abscissa (or the abscissa itself, when $\bm{a_0 = 0}$).
\end{minipage}}
& \multicolumn{1}{c}{\begin{minipage}{6.6em}
\scriptsize
\begin{center}
\vskip-3.15cm
$\bm{a_2 = 0}$ and $\bm{a_1 > 0}$
\end{center}
\vskip-.1cm
The right-hand side $q_2(x)$ of (\ref{split}) is a straight line with negative slope intersecting the abscissa at $- a_0/a_1$ and the $y$-axis at $-a_0$.
\end{minipage}}
\\
\end{tabular}
\end{center}

\begin{center}
\begin{tabular}{cccc}
\hskip7.5mm \includegraphics[width=27mm]{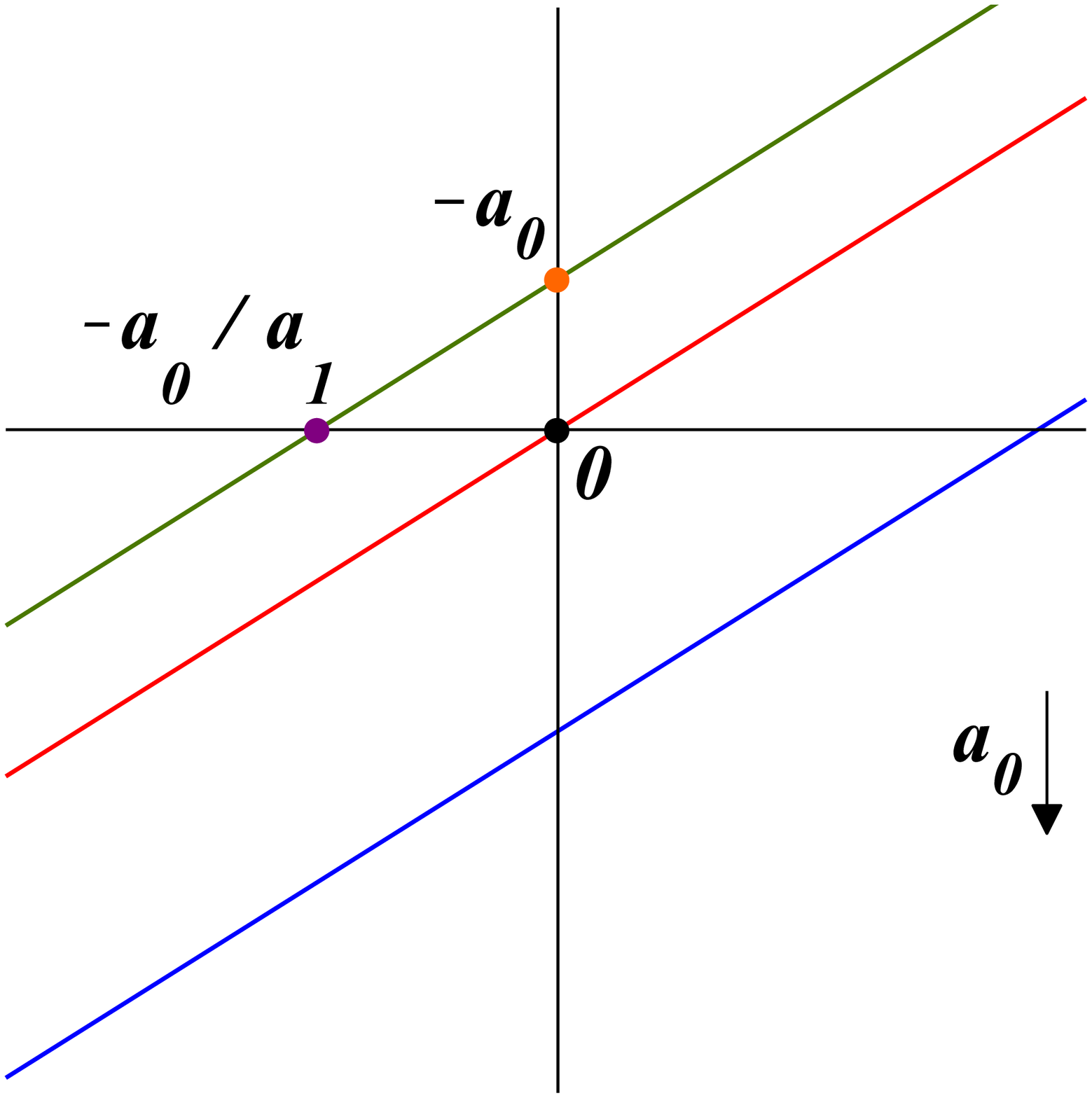} \hskip5mm & \hskip5mm \includegraphics[width=27mm]{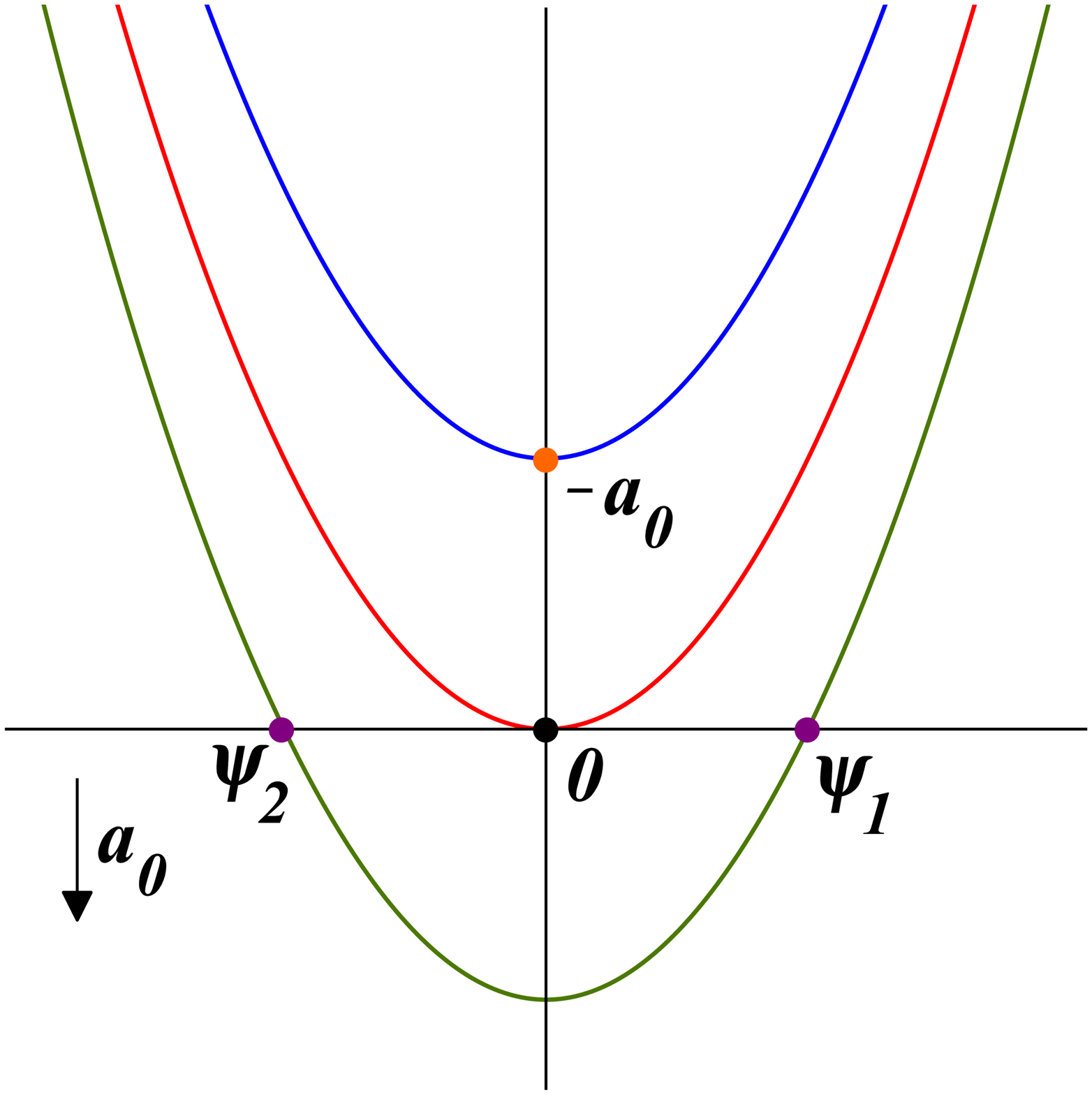} \hskip5mm & \hskip5mm \includegraphics[width=27mm]{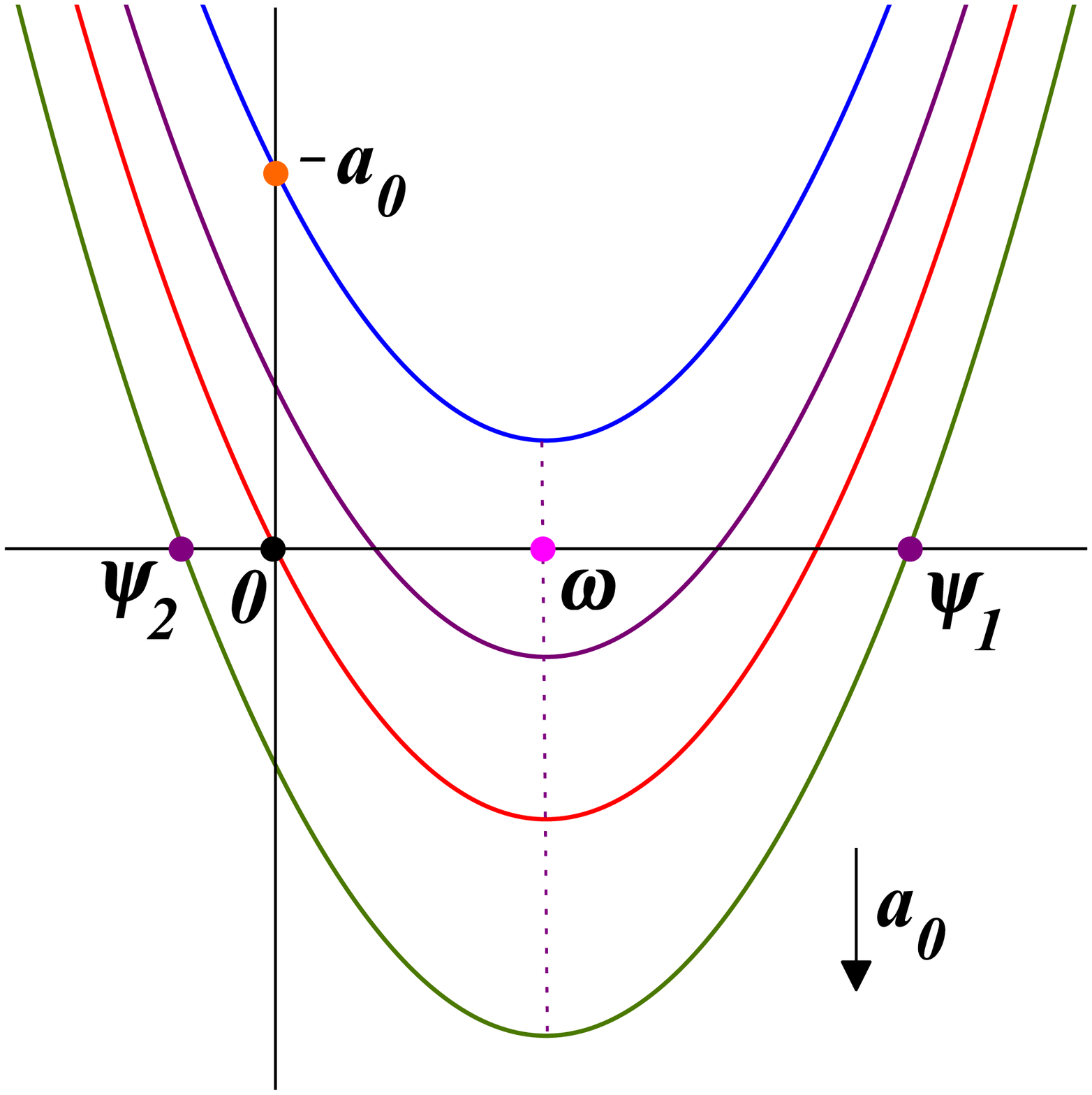} \hskip5mm & \hskip5mm \includegraphics[width=27mm]{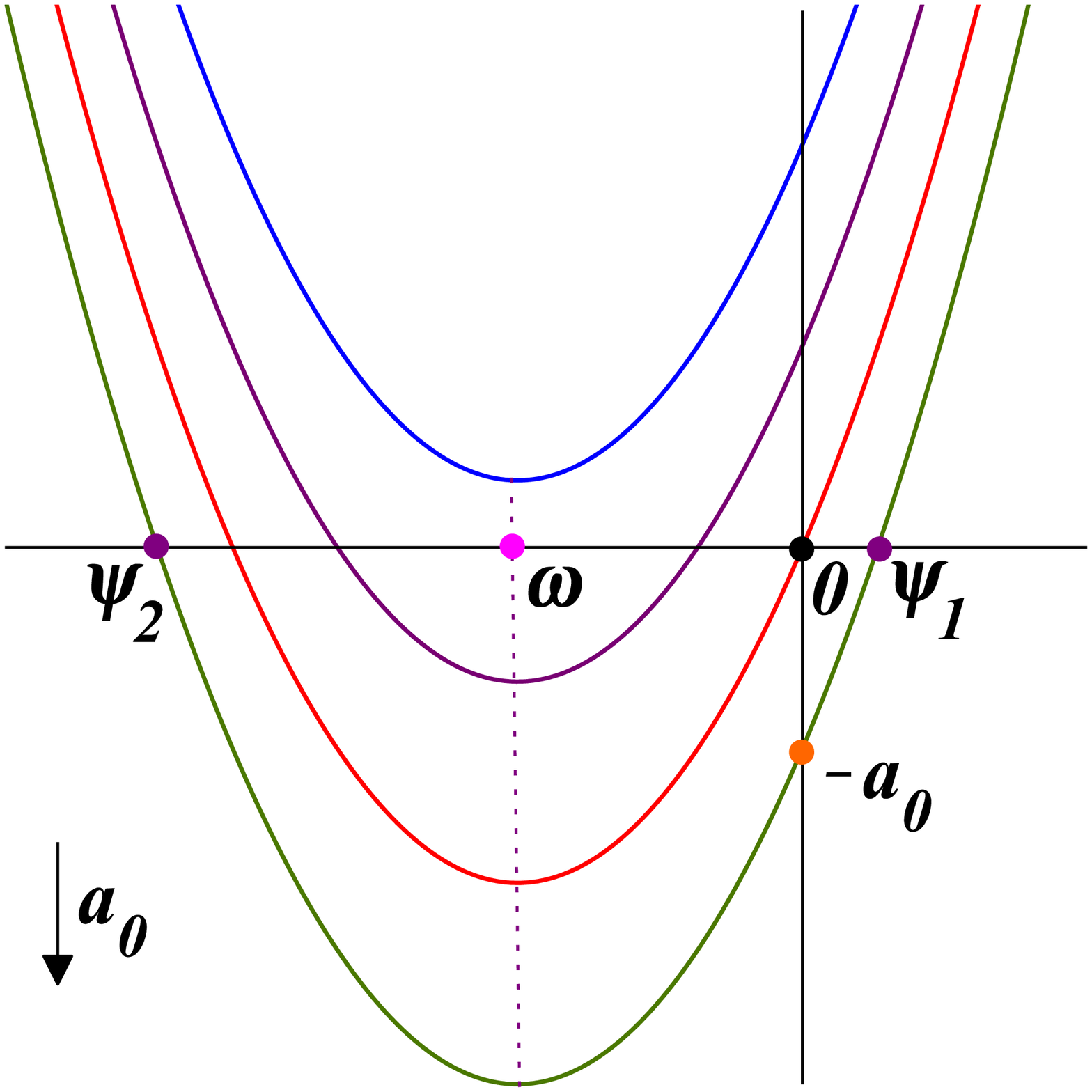} \\
{\scriptsize {\bf Figure 3f}} &  {\scriptsize {\bf Figure 3g}} & {\scriptsize {\bf Figure 3h}} &  {\scriptsize {\bf Figure 3i}} \\
& & & \\
\multicolumn{1}{c}{\begin{minipage}{6.6em}
\scriptsize
\begin{center}
\vskip-3.6cm
$\bm{a_2 = 0}$ and $\bm{a_1 < 0}$
\end{center}
\vskip-.1cm
The right-hand side $q_2(x)$ of (\ref{split}) is a straight line with positive slope intersecting the abscissa at $- a_0/a_1$ and the $y$-axis at $-a_0$.
\end{minipage}}
& \multicolumn{1}{c}{\begin{minipage}{6.6em}
\scriptsize
\begin{center}
\vskip-2.65cm
$\bm{a_2 < 0}$ and $\bm{a_1 = 0}$
\end{center}
\vskip-.1cm
If $\bm{a_0 \ge 0}$, the quadratic $q_2(x)$ has roots $\psi_{1,2} = \pm \sqrt{-a_0/a_2}$. If $\bm{a_0 < 0}$, the quadratic $q_2(x)$ has no real roots. In either case, there is a local minimum at $\omega = 0$.
\end{minipage}}
& \multicolumn{1}{c}{\begin{minipage}{6.6em}
\scriptsize
\begin{center}
\vskip-0.4cm
$\bm{a_2 < 0}$ and $\bm{a_1 > 0}$
\end{center}
\vskip-0.1cm
If $\bm{a_0 \ge a_1^2/(4 a_2)}$, the quadratic $q_2(x)$ has roots $\psi_{1,2} = - a_1/(2 a_2) \pm 1/(2 a_2)$ $\sqrt{a_1^2 - 4 a_0 a_2}$ ($\psi_1 > 0$, while $\psi_2$ could be negative, zero, or positive; $|\psi_2| < \psi_1$). If $\bm{a_0 < a_1^2/(4 a_2)}$, the quadratic $q_2(x)$ has no real roots. In either case, there is a local minimum at $\omega = -a_1/(2a_2) > 0$.
\end{minipage}}
& \multicolumn{1}{c}{\begin{minipage}{6.6em}
\scriptsize
\begin{center}
\vskip-0.4cm
$\bm{a_2 < 0}$ and $\bm{a_1 < 0}$
\end{center}
\vskip-0.1cm
If $\bm{a_0 \ge a_1^2/(4 a_2)}$, the quadratic $q_2(x)$ has roots $\psi_{1,2} = - a_1/(2 a_2) \pm 1/(2 a_2)$ $\sqrt{a_1^2 - 4 a_0 a_2}$ ($\psi_2 < 0$, while $\psi_1$ could be negative, zero, or positive; $|\psi_2| > |\psi_1|$). If $\bm{a_0 < a_1^2/(4 a_2)}$, the quadratic $q_2(x)$ has no real roots. In either case, there is a local minimum at $\omega = -a_1/(2a_2) < 0$.
\end{minipage}}
\\
\end{tabular}
\end{center}

\begin{center}
\begin{tabular}{cc}
\includegraphics[width=67mm]{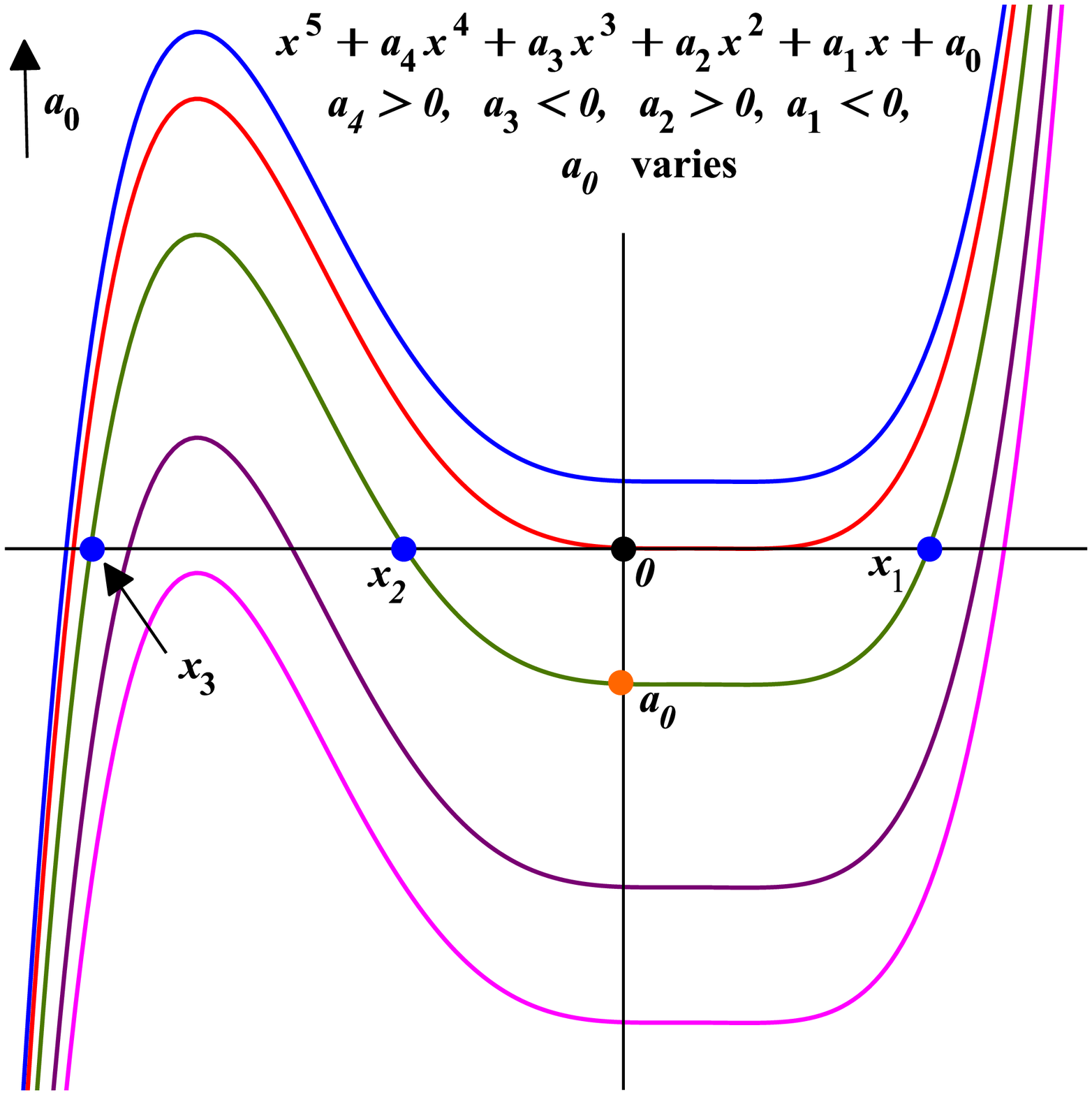} & \includegraphics[width=67mm]{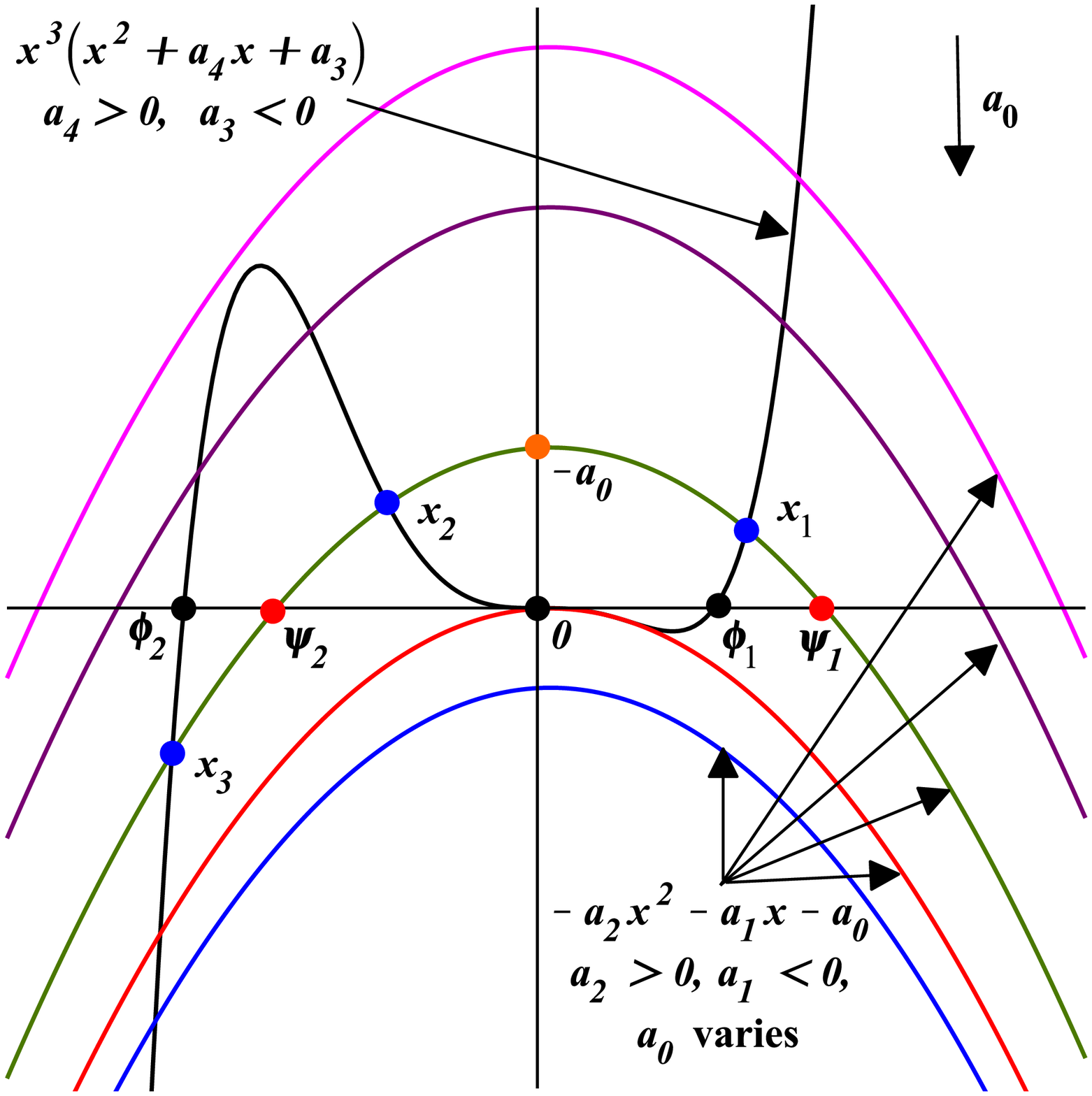} \\
& \\
{\scriptsize {\bf Figure 4a}} &  {\scriptsize {\bf Figure 4b}} \\
& \\
\multicolumn{1}{c}{\begin{minipage}{18em}
\scriptsize
\vskip-.35cm
As first example, consider a quintic with $a_4 > 0, \,\, a_3 < 0, \,\, a_2 > 0,$ and $a_1 < 0$. Shown above is the quintic $Q^{(5)}_1(x) = x^5 + x^4 - 2 x^3 + (5/6) x^2 - (1/8)x + a_0$. As the coefficient of the quadratic term, $a_2 = 5/6$, is between the roots $c_{1,2}$ of the third resolvent quadratic equation (\ref{tse}), namely, $a_2$ is between $c_2 = -3.71$ and $c_1 = 0.99$, the quintic can have zero, two, or four stationary points and hence, it can have one, or three, or five real roots. The exact number of stationary points is determined, see [4] for details, by the coefficient of the linear term, $a_1 = -1/8$. A very similar picture, impossible to discern visually from the one above, is the graph of the quintic $Q^{(5)}_2(x) = x^5 + x^4 - 2 x^3 + 3 x^2 - (1/8)x + a_0$. This quintic differs from $Q^{(5)}_1(x)$ only in the coefficient of its quadratic term --- now one has $a_2 = 3$ and such value of $a_2$ is not between $c_2 = -3.71$ and $c_1 = 0.99$. Hence, $Q^{(5)}_2(x)$ can have either zero or two stationary points and hence, it can have either one or three real roots. Again, the number of stationary points is determined [4] by the coefficient of the linear term. For both quintics, the exact number of real roots depends on the coefficient $a_0$, which is taken as a varying parameter, so that the quintics can ``slide" vertically.
\end{minipage}}
& \multicolumn{1}{c}{\begin{minipage}{18em}
\scriptsize
\vskip-2.35cm
The two ``components" of the quintic $Q^{(5)}_1(x)$ [or $Q^{(5)}_2(x)$] from Figure 4a, namely, the ``sub-quintic" $x^3 (x^2 + x - 2)$ (same for both of them) and the quadratic $- a_2 x^2 + (1/8) x - a_0$, where $a_2 = 5/6$ for $Q^{(5)}_1(x)$ [or $a_2 = 3$ for $Q^{(5)}_1(x)$] --- again, impossible to discern visually one from the other. Varying the quintic coefficient $a_0$ results in ``sliding" the quadratics $- a_2 x^2 + (1/8) x - a_0$ vertically. If $a_2 = 5/6$, the sliding quadratic $- (5/6) x^2 + (1/8) x - a_0$ can be tangent to the ``sub-quintic" $x^3 (x^2 + x - 2)$ at 4 points. Thus, the quintic $Q^{(5)}_1(x)$ can have 1, or 3, or 5 real roots. If $a_2 = 3$, the sliding quadratic $-3 x^2 + (1/8) x - a_0$ can be tangent to the ``sub-quintic" $x^3 (x^2 + x - 2)$ at 2 points only. The quintic $Q^{(5)}_2(x)$ can therefore have either 1 or 3 real roots only. The points at which the quadratic $- a_2 x^2 + (1/8) x - a_0$ can be tangent to the ``sub-quintic" $x^3 (x^2 + x - 2)$ are the roots $\xi_i$ of the auxiliary quartic equation (\ref{quartic}).
\end{minipage}}
\\
\end{tabular}
\end{center}
This is the full analysis of $Q^{(5)}_1$:
\begin{itemize}
\item [(i)] When $a_0$ is sufficiently large [$a_0 = 1$  is sufficiently large for $Q^{(5)}_1(x) = x^5 + x^4 - 2 x^3 + (5/6) x^2 - (1/8) x + a_0$], the only real root of the resulting quintic $x^5 + x^4 - 2 x^3 + (5/6) x^2 - (1/8) x + 1$ is the negative $x = -2.16$ --- as the quadratic $-(5/6) x^2 + (1/8) x - 1$ does not ``reach" the ``sub-quintic" $x^3(x^2 + x - 2)$ in the fourth quadrant --- see the lowermost quadratic on Figure 4b. The quadratic $-(5/6) x^2 + (1/8) x - 1$ does not have real roots. \\
    One can find a sufficient condition for not having positive real roots: if the local minimum of the ``sub-quintic" (at  $x > 0$) is greater than the absolute maximum of the parabola, then there can be no positive roots of the quintic. In this particular case the local minimum of the ``sub-quintic" $x^3(x^2 + x - 2)$ is $-0.29$, while the absolute maximum of the quadratic $-(5/6) x^2 + (1/8) x - 1$ is $-0.9953$.
    
\item [(ii)] ``Sliding" the quadratic $- (5/6) x^2 + (1/8) x - a_0$ up by decreasing $a_0$ will eventually allow the parabola to ``reach" the ``sub-quintic" in the fourth quadrant. Taking $a_0 = 1/100$, yields the quintic $x^5 + x^4 - 2 x^3 + (5/6) x^2 - (1/8) x + 1/100$. The real roots of this quintic are three --- one negative and two positive: $x_3 = -2.13$, $x_2 = 0.44$, and $x_1 = 0.51$.
    
\item [(iii)] Bringing the quadratic $-(5/6) x^2 + (1/8) x - a_0$ further up results in the appearance of two additional positive roots. Taking $a_0 = 6/1000$, gives the quintic $x^5 + x^4 - 2 x^3 + (5/6) x^2 - (1/8) x + 6/1000$ which has five real roots --- one negative and four positive: $x_5 = -2.13$, $x_4 = 0.10$, $x_3 = 0.17$, $x_2 = 0.30$, and $x_1 = 0.56$. The fact that the roots of the quadratic $-(5/6) x^2 + (1/8) x - 6/1000$ are not real is insignificant. That is, one can still have a quintic with five real roots (one negative and four positive) when the quadratic has two real roots. For example, for the quintic $x^5 + x^4 - 2 x^3 + (5/6) x^2 - (5/48) x + 1/10000$, the roots are $x_5 = -2.13$, $x_4 = 0.0010$, $x_3 = 0.26$, $x_2 = 0.37$, and $x_1 = 0.50$. The resulting quadratic $-(5/6) x^2 + (5/48) x - 1/10000$ has two real roots: and $\psi_2 = 0.0010$ and $\psi_1 = 0.12$. \\
    For as long as the bigger root $\psi_1$ of the parabola is smaller than the bigger root $\phi_1$ of the ``sub-quintic", the number of positive roots (in addition to the unique negative root) is 0, or 2, or 4. \\
    Note that with the decrease of the free term of the auxiliary quartic equation (\ref{quartic}), that is, with the decrease of the linear term $a_1$ of the quintic, two of the stationary points of the quintic will eventually disappear, i.e. there will be only two real roots of the auxiliary quartic equation (\ref{quartic}) --- see the analysis on Figure 1.13 in [4]. Hence, it will not be possible to have a quintic with five real roots (with the decrease of the quintic's linear term $a_1$, two of the roots of the auxiliary quartic equation will coalesce into a double root and then become complex). As example in this vein, alter the above quintic $x^5 + x^4 - 2 x^3 + (5/6) x^2 - (1/8) x + 6/1000$ (which has five real roots), into the quintic $x^5 + x^4 - 2 x^3 + (5/6) x^2 - 3 x + 6/1000$. Two of the positive roots of the quintic become complex and the real roots are now only three: the negative $x_3 = -2.29$ and the two positive $x_2 = 0.002$ and $x_1 = 1.32$. \\
    One can also note from Figure 1.13 in [4] that if the free term of the auxiliary quartic equation (\ref{quartic}) is negative, as it is in the current case, it is not possible to have a situation with no stationary points of the quintic.
    
    \vskip.5cm
    \n
    One can summarize that for sub-cases (i), (ii), and (iii), when the bigger root $\psi_1$ of the parabola is smaller than the bigger root $\phi_1$ of the ``sub-quintic", there is a unique negative root smaller than the smaller root $\phi_2$ of the ``sub-quintic" and there also are 0, or 2, or 4 (subject to the conditions discussed) positive roots smaller than the bigger root $\phi_1$ of the ``sub-quintic". \\
    One can also determine a lower bound for the negative root of all quintics in (i), (ii), and (iii). Firstly, one has to replace $x$ by $-x$ in the quintic and then seek the upper root bound. For example, for the quintic $x^5 + x^4 - 2 x^3 + (5/6) x^2 - (1/8) x + 1$, changing $x$ to $-x$ yields $-(x^5 - x^4 - 2 x^3 - (5/6) x^2 - (1/8) x - 1)$ whose roots are the same as those of $x^5 - x^4 - 2 x^3 - (5/6) x^2 - (1/8) x - 1$. An upper bound, as proposed in [5], is the bigger of 1 and the sum of the absolute values of all negative coefficients: $1 + 2 + 5/6 + 1/8 + 1 = 119/24 \approx 4.96$. [As $x^5 - x^4 - 2 x^3 - (5/6) x^2 - (1/8) x - 1$ is a Cauchy polynomial, see [1], i.e. all of its coefficients are negative except the leading positive coefficient,  in this case the bound [5] coincides with the Cauchy bound.] Hence, the lower bound for the roots of the quintic $x^5 + x^4 - 2 x^3 + (5/6) x^2 - (1/8) x + 1$ is $-4.96$. One can also use, for example, the Kurosh bound [6]: $1 + \sqrt[k]{B/a_0}$, where $a_k$ (with $k \ge 1$) is, in Kurosh notation, the first negative coefficient (in the sequence of decreasing powers of $x$; the leading term having positive coefficient $a_0$) and $B$ is the biggest absolute value of all negative coefficients. For the quintic $x^5 - x^4 - 2 x^3 - (5/6) x^2 - (1/8) x - 1$, one has $k = 1$ and $B = 2$ and, hence, the following upper root bound: $1 + \sqrt[1]{2} = 1 + 2 = 3$. The lower root bound for the quintic $x^5 + x^4 - 2 x^3 + (5/6) x^2 - (1/8) x + 1$ is therefore $-3$.
    
\item [(iv)] If $a_0 > 0$ and the bigger root of the parabola, $\psi_1$, is greater than the bigger root $\phi_1$ of the ``sub-quintic" (see Figure 4b), then this is a sufficient condition for not being possible to have four positive roots. That is, the number of positive roots of the quintic (in addition to the single negative root) could be either zero or two. As example of the former, consider the quintic $x^5 + x^4 - 2 x^3 + (5/6) x^2 - 3x + 5/2$. For it, one has $\psi_2 = 1.31$ and $\psi_1 = 2.29$, while $\phi_1 = 1.00$. The only real root is $x = -2.33$. As example for the latter, consider the quintic $x^5 + x^4 - 2 x^3 + (5/6) x^2 - 3x + 1$ comes with $\psi_2 = 1.37$ and $\psi_1 = 3.23$, while $\phi_1 = 1.00$ again. The roots are: $x_3 = -2.30$, $x_2 = 0.35$, and $x_1 = 1.24$
    
\item [(v)] Taking the ``separatrix" quadratic (corresponding to $a_0 = 0$), ``slides" one of the positive roots of the quintic to the origin. In result, one can have either a quintic with one negative root, a zero root and a positive root: $x^5 + x^4 - 2 x^3 + (5/6) x^2 - (1/8) x$ (with roots: $x_3 = -2.13$, $x_2 = 0$, and $x_1 = 0.60$), or a quintic with one negative root, a zero root and three positive roots: $x^5 + x^4 - 2 x^3 + (5/6) x^2 - (5/48) x$ (with roots: $x_5 = -2.13$, $x_4 = 0$, $x_3 = 0.27$, $x_2 = 0.36$, and $x_1 = 0.50$). The unique negative root is smaller than the smaller root $\phi_2$ of the ``sub-quintic" (one needs to find a lower root bound in order to determine its isolation interval). The one or three positive roots are smaller than the bigger root $\phi_1$ of the ``sub-quintic".

\item [(vi)] Taking next a negative value of $a_0$, ``slides" the zero root into the negative part of the abscissa and hence, in result, one can have either a quintic with two negative roots and one positive root: $x^5 + x^4 - 2 x^3 + (5/6) x^2 - (1/8) x - 1/2$ (with roots: $x_3 = -2.12$, $x_2 = -0.48$, and $x_1 = 0.95$), or a quintic with two negative roots and three positive roots: $x^5 + x^4 - 2 x^3 + (5/6) x^2 - (1/10) x - 1/1000$ (with roots: $x_5 = -2.13$, $x_4 = -0.009$, $x_3 = 0.26$, $x_2 = 0.42$, and $x_1 = 0.46$). \\
    For as long as the bigger root $\psi_1$ of the parabola is smaller than the bigger root $\phi_1$ of the ``sub-quintic", the only positive root or the three positive roots are greater than the bigger root $\psi_1$ of the parabola and smaller than the bigger root $\phi_1$ of the ``sub-quintic". \\
    One of the two negative roots of the quintic is smaller than the smaller root $\phi_2$ of the ``sub-quintic" (one needs to find a lower root bound in order to determine its isolation interval) and the other negative root of the quintic is greater than the smaller root $\psi_2$ of the parabola. \\
    As in sub-case (iii) above, with the decrease of the free term of the auxiliary quartic equation (\ref{quartic}), that is, with the decrease of the linear term $a_1$ of the quintic, two of the stationary points of the quintic will eventually disappear, i.e. there will be only two real roots of the auxiliary quartic equation (\ref{quartic}) --- see, again, the analysis on Figure 1.13 in [4]. Hence, it will not be possible to have a quintic with five real roots. Again, one can also note from Figure 1.13 in [4] that if the free term of the auxiliary quartic equation (\ref{quartic}) is negative, it is not possible to have a situation with no stationary points of the quintic. As example, alter the above quintic $x^5 + x^4 - 2 x^3 + (5/6) x^2 - (1/10) x - 1/1000$ (which has five real roots), into the quintic $x^5 + x^4 - 2 x^3 + (5/6) x^2 - (1/8) x - 1/1000$. Two of the positive roots of the quintic become complex and the real roots are now only three: the negative $x_3 = -2.14$ and the two positive $x_2 = 0.008$ and $x_1 = 0.60  $. \\
    As in (iv) above, if the bigger root of the parabola, $\psi_1$, is greater than the bigger root $\phi_1$ of the ``sub-quintic", this is a sufficient condition to eliminate the possibility of having three positive roots --- one can have only one positive root in such case. As example, the quintic  $x^5 + x^4 - 2 x^3 + (5/6) x^2 - (1/8) x - 1$ comes with $\psi_2 = -1.02$ and $\psi_1 = 1.17 > \phi_1 = 1.00$ and has two negative roots: $x_3 = -2.10$ and $x_2 = -0.64$ and one positive root: $x_1 = 1.06$.
   
\item [(vii)] Decreasing $a_0$ further, ``slides" the parabola upwards. While the smaller root of the parabola, $\psi_2$, is greater than the smaller root $\phi_2$ of the ``sub-quintic" (see Figure 4b), then, in addition to either three or one positive roots, there will necessarily be two negative roots, one of which is greater than $\psi_2$ and the other --- smaller than $\phi_2$ (one needs a lower root bound for the determination of the isolation interval of the latter). For the already considered quintic $x^5 + x^4 - 2 x^3 + (5/6) x^2 - (1/10) x - 1/1000$, one has: $\psi_1 = 0.13$, $\psi_2 = -0.00928$, $\phi_1 = 1.00$, and $\phi_2 = -2.00$. Hence, one of its negative roots (in addition to the three positive roots $x_3 = 0.26$, $x_2 = 0.42$, and $x_1 = 0.46$) is greater than $\psi_2 = -0.00928$, namely: $x_4 = -0.00926$, and the other negative root is smaller than $\phi_2 = -2.00$, that is, $x_5 = -2.13$.
    
\item [(viii)] Bringing this parabola further up, results in two of the three positive roots of the quintic coalescing into a double root and then becoming complex. The quintic, in addition to its two negative roots (should the smaller root of the parabola, $\psi_2$, be still  greater than the smaller root $\phi_2$ of the ``sub-quintic"), will now only have one positive root --- between min$\{\psi_1, \phi_1\}$ and max$\{\psi_1, \phi_1\}$.  As example, take $a_0 = -1/2$ and ``slide" the quintic of sub-case (vii) up to $x^5 + x^4 - 2 x^3 + (5/6) x^2 - (1/10) x - 1/2$. One now has: $\psi_1 = 0.84$, $\psi_2 = -0.72$, $\phi_1 = 1.00$, and $\phi_2 = -2.00$.  The unique positive root is $x_1 = 0.94$ --- between $\psi_1 = 0.84$ and $\phi_1 = 1.00$. One of the two negative roots, namely $x_2 = -0.48$, is between $\psi_2 = -0.72$ and $0$. The other negative root, $x_3 = -2.12$, is smaller than $\phi_2 = -2.00$. Again, a lower bound for the smallest roots of quintics, such as those in sub-cases (vii) and (viii), is needed in order to determine the root isolation intervals of their smallest roots.
    
\item [(ix)] Once $a_0$ is small enough, the smaller root $\psi_2$ of the parabola will become smaller than the smaller root of the ``sub-quintic" --- see the uppermost parabola on Figure 4b. Then, in addition to the only positive root, there will be either two negative roots greater than the smaller root $\phi_2$ of the ``sub-quintic" or no negative roots at all. \\
    As example of the former, consider the quintic $x^5 + x^4 - 2 x^3 + (5/6) x^2 - (1/8) x - 4$. One has: $\psi_1 = 2.27$, $\psi_2 = -2.12$, $\phi_1 = 1.00$, and $\phi_2 = -2.00$. The positive root is $x_1 = 1.34$ and the negative roots are $x_2 = -1.15$ and $x_3 = -1.98$ --- between $\phi_2 = -2.00$ and zero. \\
    The quintic $x^5 + x^4 - 2 x^3 + (5/6) x^2 - (1/8) x - 8$ has $\psi_1 = 3.17$, $\psi_2 = -3.02$, $\phi_1 = 1.00$, and $\phi_2 = -2.00$. The positive root is $x_1 = 1.34$ and there are no negative roots. \\
    One can find a sufficient condition for having two negative real roots: if the local maximum of the ``sub-quintic" (at  $x < 0$) is greater than the absolute maximum of the parabola, then there will be two negative roots. In the particular case of the quintic $x^5 + x^4 - 2 x^3 + (5/6) x^2 - (1/10) x - 1/2$ from sub-case (viii), the local maximum of the ``sub-quintic" $x^3(x^2 + x - 2)$ is $4.27$, while the absolute maximum of the quadratic $-(5/6) x^2 + (1/10) x + 1/2$ is $0.503$. Hence, the quintic has two negative roots (in addition to the unique positive root).

\end{itemize} 

\vskip.5cm
\n
Consider next the quintic $Q^{(5)}_2(x) = x^5 + x^4 - 2 x^3 + 3 x^2 - (1/8)x + a_0$. \\ 
The graph of $Q^{(5)}_2(x)$ is similar to the one of $Q^{(5)}_1(x)$ on Figure 4a --- it is impossible to visually discern one from the other. The quintic $Q^{(5)}_2(x)$ differs from $Q^{(5)}_1(x)$ only in the coefficient of its quadratic term --- now one has $a_2 = 3$ and such value of $a_2$ is not between $c_2 = -3.71$ and $c_1 = 0.99$.  Hence, the quintic $Q^{(5)}_2(x)$ can have either zero or two stationary points and hence, it can have either one or three real roots. \\
The full analysis of the quintic $Q^{(5)}_2$ is as follows:

\begin{itemize}
\item [(I)] When $a_0$ is sufficiently large [$a_0 = 1/2$  is sufficiently large for $Q^{(5)}_2(x) = x^5 + x^4 - 2 x^3 + 3 x^2 - (1/8) x + a_0$], the only real root of the resulting quintic $x^5 + x^4 - 2 x^3 + 3 x^2 - (1/8) x + 1/2$ is the negative $x = -2.39$ --- as the quadratic $-3 x^2 + (1/8) x - 3$ does not ``reach" the ``sub-quintic" $x^3(x^2 + x - 2)$ --- see the lowermost quadratic on Figure 4b. The quadratic $-(5/6) x^2 + (1/8) x - 1/2$ does not have real roots. \\
    As in sub-case (i) for the quintic $Q^{(5)}_1(x)$, one has the same sufficient condition for not having positive real roots: if the local minimum of the ``sub-quintic" (at  $x > 0$) is greater than the absolute maximum of the parabola, then there can be no positive roots of the quintic. In this particular case the local minimum of the ``sub-quintic" $x^3(x^2 + x - 2)$ is $-0.29$, while the absolute maximum of the quadratic $-3 x^2 + (1/8) x - 1/2$ is $-0.4987$.

\item [(II)] ``Sliding" the quadratic $- 3 x^2 + (1/8) x - a_0$ up by decreasing $a_0$ will eventually allow the parabola to ``reach" the ``sub-quintic", yielding two positive roots in addition to the unique negative root. The intersection of the parabola with the ``sub-quintic" can happen in the fourth quadrant only, or in the first quadrant only, or both in the fourth and first quadrant: \\
    Taking $a_0 = 1/760$, yields the quintic $x^5 + x^4 - 2 x^3 + 3 x^2 - (1/8) x + 1/760$. The parabola $-3 x^2 + (1/8) x - 1/760$ has no real roots, but it ``reaches" the ``sub-quintic" in the fourth quadrant. The real roots of the quintic are three --- one negative and two positive: $x_3 = -2.38$, $x_2 = 0.0020$, and $x_1 = 0.0023$. \\
    Taking $a_0 = 1/10000$, yields the quintic $x^5 + x^4 - 2 x^3 + 3 x^2 - (1/8) x + 1/10000$. The parabola $-3 x^2 + (1/8) x - 1/10000$ has two real roots: $0.0008$ and $0.0408$. This parabola also ``reaches" the ``sub-quintic" in the fourth quadrant. The real roots of the quintic are three --- one negative and two positive: $x_3 = -2.38$, $x_2 = 0.0008$, and $x_1 = 0.0420$. \\
    Taking $a_0 = 25/2$ and $a_1 = 15$, yields the quintic $x^5 + x^4 - 2 x^3 + 3 x^2 - 15 x + 25/2$. The parabola $-3 x^2 + 15 x - 25/2$ has two real roots: $1.06$ and $3.94$. The parabola ``reaches" the ``sub-quintic" in the first quadrant only. The real roots of the quintic are three --- one negative and two positive: $x_3 = -2.87$, $x_2 = 1.12$, and $x_1 = 1.24$. \\
    Taking $a_0 = 10$ and $a_1 = 15$, yields the quintic $x^5 + x^4 - 2 x^3 + 3 x^2 - 15 x + 10$. The parabola $-3 x^2 + 15 x - 10$ has two real roots: $0.79$ and $4.21$. The parabola ``reaches" the ``sub-quintic" in both the fourth quadrant and in the first quadrant. The real roots of the quintic are three --- one negative and two positive: $x_3 = -2.86$, $x_2 = 0.76$, and $x_1 = 1.49$. 
    
\vskip.5cm
\n
One can summarize that for cases (I) and (II) (positive $a_0$), when the discriminant of the parabola is negative (the parabola having no real roots), the quintic has a unique negative root smaller than the smaller root $\phi_2$ of the ``sub-quintic" (a lower root bound is needed for its isolation interval) and there are either zero or two positive roots smaller than the bigger root $\phi_1$ of the ``sub-quintic". If the local minimum of the ``sub-quintic" (at  $x > 0$) is greater than the absolute maximum of the parabola, then there can be no positive roots of the quintic. \\
If the discriminant of the parabola is non-negative (the parabola will then have two real positive roots) and the smaller root $\psi_2$ of the parabola is smaller than the bigger root $\phi_1$ of the ``sub-quintic", then the quintic will have a unique negative root smaller than the smaller root $\phi_2$ of the ``sub-quintic" (a lower root bound is needed for its isolation interval) and there will be two positive roots: one of them smaller than the smaller root $\psi_2$ of the parabola and the other --- between the smaller of the bigger root $\psi_1$ of the parabola and the bigger root $\phi_1$ of the ``sub-quintic" and the bigger of the two, i.e. between min$\{\psi_1, \phi_1\}$ and max$\{\psi_1, \phi_1\}$. \\
If the parabola has a non-negative discriminant and the smaller root $\psi_2$ of the parabola is greater than the bigger root $\phi_1$ of the ``sub-quintic", then the quintic will have a unique negative root smaller than the smaller root $\phi_2$ of the ``sub-quintic" (with lower root bound needed for the determination of its isolation interval) and there will be either zero or two positive roots and these would be either between the two roots $\psi_2$ and $\psi_1$ of the parabola or would be smaller than the bigger root $\phi_1$ of the ``sub-quintic". \\

\item [(III)]  Taking the ``separatrix" quadratic (corresponding to $a_0 = 0$), ``slides" one of the two positive roots of the quintic to the origin. In result, one will have a quintic with one negative root, a zero root and one positive root. The negative root is smaller than the smaller root $\phi_2$ of the ``sub-quintic" (a lower root bound is needed for its isolation interval) and the positive root is between the smaller of the bigger root $\psi_1$ of the parabola and the bigger root $\phi_1$ of the ``sub-quintic" and the bigger of the two, i.e. between min$\{\psi_1, \phi_1\}$ and max$\{\psi_1, \phi_1\}$. \\
    For example, the quintic $x^5 + x^4 - 2 x^3 + 3 x^2 - (1/8) x$ comes with $\psi_2 = 0$, $\psi_1 = 0.04167$, $\phi_2 = -2$, and $\phi_1 = 1$. It should have a negative root smaller than $\phi_2 = -2$, a zero root, and a positive root between min$\{\psi_1, \phi_1\}$ and max$\{\psi_1, \phi_1\}$. Indeed, the roots of the quintic are $x_3 = -2.38$, $x_2 = 0$, and $x_1 = 0.04286$ --- exactly as predicted. 

\item [(IV)] Taking next a negative value of $a_0$, ``slides" the zero root of the separatrix quintic into the negative part of the abscissa and, for as long as the smaller root $\psi_2$ of the parabola is greater than or equal to the smaller root $\phi_2$ of the ``sub-quintic", the quintic will have two negative roots and one positive root. The smaller negative roots is smaller than the smaller root $\phi_2$ of the ``sub-quintic" (a lower root bound is needed for its isolation interval), the bigger negative root is greater than the smaller root $\psi_2$ of the parabola. The positive root is between the smaller of the bigger root $\psi_1$ of the parabola and the bigger root $\phi_1$ of the ``sub-quintic" and the bigger of the two. \\
    As example, consider the quintic $x^5 + x^4 - 2 x^3 + 3 x^2 - (1/8) x - 1$. It has:  $\psi_2 = -0.56$, $\psi_1 = 0.60$, $\phi_2 = -2$, and $\phi_1 = 1$. Clearly, $\phi_2 = -2 < \psi_2 = -0.56$. Hence, there should be a root smaller than $-2.00$, a root between $-0.56$ and 0, and a root between $0.60$ and $1$. Indeed, the roots are: $x_3 = -2.36$, $x_2 = -0.48$, and $x_1 = 0.67$. \\
    
\item [(V)] Sliding the parabola further up (by decreasing $a_0$) will eventually result in the smaller root $\psi_2$ of the parabola becoming smaller than the smaller root $\phi_2$ of the ``sub-quintic". The resulting quintic will then have one positive root between the smaller of the bigger root $\psi_1$ of the parabola and the bigger root $\phi_1$ of the ``sub-quintic" and the bigger of the two, i.e. between min$\{\psi_1, \phi_1\}$ and max$\{\psi_1, \phi_1\}$, and either no other roots or two negative roots --- both greater than the smaller root  $\phi_2$ of the ``sub-quintic". \\
    As example of the latter, consider the quintic $x^5 + x^4 - 2 x^3 + 3 x^2 - (1/8) x - 13$. It has $\psi_2 = -2.06$, $\psi_1 = 2.10$, $\phi_2 = -2$, and $\phi_1 = 1$. Clearly, $\psi_2 = -2.06 < \phi_2 = -2$. The quintic has one positive real root ($x_1 = 1.52$) and two negative roots greater than $\phi_2 = -2$. These are $x_2 = -1.73$ and $x_3 = - 1.91$. \\
    As example of the former, consider the quintic $x^5 + x^4 - 2 x^3 + 3 x^2 - (1/8) x - 14$. It has $\psi_2 = -2.14$, $\psi_1 = 2.18$, $\phi_2 = -2$, and $\phi_1 = 1$. Clearly, $\psi_2 = -2.14 < \phi_2 = -2$. The quintic has one real root only --- the positive $x_1 = 1.54$ --- between $\phi_1$ and $\psi_1$.
\end{itemize}

\vskip.5cm
\subparagraph{\hskip-.6cm 5 \hskip0.2cm Discussion \vskip.5cm}
\hskip-1cm  It is possible to significantly reduce the number of cases for the general quintic $Q(x) = x^5 + a_4 x^4 + a_3 x^3 + a_2 x^2 + a_1 x + a_0$ by depressing it, that is, by eliminating the quartic term with a suitable change of the coordinate system: $x \to x - a_4/5$. The resulting depressed quintic is $x^5 + p x^3 + q x^2 + r x + s$, where:
\b
\label{pe}
p & = & -\frac{2}{5} a_4^2 + a_3, \\
q & = &  \frac{4}{25} a_4^3 - \frac{3}{5} a_3 a_4 + a_2, \\
r & = & -\frac{3}{125} a_4^4 + \frac{3}{25} a_3 a_4^2 - \frac{2}{5} a_2 a_4 + a_1 , \\
\label{es}
s & = & \frac{4}{3 125} a_4^5 - \frac{1}{125} a_3 a_4^3 + \frac{1}{25} a_2 a_4^2 - \frac{1}{5} a_1 a_4 + a_0 .
\e
The depressed ``sub-quintic" involved now is $x^5 + p x^3$ and it has, instead of 15 possible cases, only three cases (two of which are qualitatively different) --- see Figures 2a, 2b, and 2c with $a_3$ replaced by $p = -(2/5) a_4^2 + a_3$. \\
The example ``sub-quintic", analyzed in the previous section, has $a_4 = 1$ and $a_3 = -2$. Hence, $p = -12/5$ and Figure 2c applies now (as well as Figure 3c for the parabola). One can immediately see that the reduction of the total number of possible cases by depressing the quintic does not provide much benefit: there are hardly any gains associated with the graph of the ``sub-quintic" becoming centrally symmetric (with respect to the origin), as opposed to the ``original" graph on Figure 2i. Separately, the coefficients of the parabola become too complex and the analysis of the roots of the parabola and the interrelation between the coefficients of the quintic becomes very tedious. \\
One can also eliminate the cubic term of the quintic by making the transformation $x \to x - \beta$, where $\beta$ is either of the roots of the quadratic equation $10 \beta^2 -4 a_4 \beta + a_3 = 0$. Then there will be only two remaining cases for the resulting depressed ``sub-quintic" --- Figures 2d and 2j. However, the coefficients of the depressed quintic without the cubic term would be even more complicated and the analysis will be even more tedious. \\
However, despite the complexity of the conditions of the complete root classification of the quintic, (\ref{de2})--(\ref{ef2}), one can eliminate all indeterminacy of the proposed method, by calculating $D_2, \,\, D_3, \,\, D_4, \,\, D_5, \,\, E_2,$ and $F_2$ for $p, \,\ q, \,\, r,$ and $s$ given by (\ref{pe})--(\ref{es}). Then the {\it exact} number of roots of the quintic within the intervals determined by the roots of the two {\it resolvent} quadratic polynomials $q_1(x) = x^2 + a_4 x + a_3$ and $q_2(x) = a_2 x^2 + a_1 x + a_0$ will be known.

\pagebreak 
\n
\subparagraph{\hskip-.6cm References \vskip.2cm}
\n
\begin{table}[!h]
\begin{tabular}{ll}
\multicolumn{1}{l}{\begin{minipage}{0.1cm} \vskip-.5cm $\!\!$[1] \end{minipage}}
& \multicolumn{1}{c}{\begin{minipage}{14.3cm} Q.I. Rahman and G. Schmeisser, {\it Analytic Theory of Polynomials}, Oxford University Press (2002); \end{minipage}} \\ 
\multicolumn{1}{l}{\begin{minipage}{0.1cm} \vskip-.5cm $\!\!$ \end{minipage}}
& \multicolumn{1}{c}{\begin{minipage}{14.3cm} V.V. Prasolov, {\it Polynomials}, Springer (2010); \end{minipage}} \\
\multicolumn{1}{l}{\begin{minipage}{0.1cm} \vskip-.5cm $\!\!$ \end{minipage}}
& \multicolumn{1}{c}{\begin{minipage}{14.3cm} R. Bruce King, {\it Beyond the Quartic Equation}, Birkh\"auser (1996); \end{minipage}} \\
\multicolumn{1}{l}{\begin{minipage}{0.1cm} \vskip-.5cm $\!\!$ \end{minipage}}
& \multicolumn{1}{c}{\begin{minipage}{14.3cm} M.M. Postnikov, {\it Foundations of Galois Theory}, Dover (2004). \end{minipage}} \\ \\
\multicolumn{1}{l}{\begin{minipage}{0.1cm} \vskip-0.7cm $\!\!$[2] \end{minipage}}
& \multicolumn{1}{c}{\begin{minipage}{14.3cm} \vskip-.2cm D.S. Arnon, {\it Geometric Reasoning with Logic
and Algebra}, Artificial Intelligence {\bf 37}, 37--60 (1988); \end{minipage}} \\  \\
\multicolumn{1}{l}{\begin{minipage}{0.1cm} $\!\!$ \end{minipage}}
& \multicolumn{1}{c}{\begin{minipage}{14.3cm} \vskip-0.4cm J.R. Johnson, {\it Algorithms for Polynomial Real Root Isolation (PhD thesis, 1991)}, in: B.F. Caviness and J.R. Johnson (Eds.), {\it Quantifier Elimination and Cylindrical Algebraic Decomposition}, Springer (1998); \end{minipage}} \\ \\
\multicolumn{1}{l}{\begin{minipage}{0.1cm} \vskip-.7cm $\!\!$ \end{minipage}}
& \multicolumn{1}{c}{\begin{minipage}{14.3cm} \vskip-.4cm  Lu Yang, H. Xiaorong, and Zh. Zeng, {\it Complete Discrimination System for Polynomials}, Science in China (Series E, Technological Sciences) {\bf 39(6)}, 628--646 (1996); \end{minipage}} \\ \\
\multicolumn{1}{l}{\begin{minipage}{0.1cm} \vskip-1.5cm $\!\!$ \end{minipage}}
& \multicolumn{1}{c}{\begin{minipage}{14.3cm} \vskip-0.4cm S. Liang and D.J. Jeffrey, {\it An Algorithm for Computing the Complete
Root Classification of a Parametric Polynomial}, in: J. Calmet, T. Ida, and D. Wang (Eds.), {\it Artificial Intelligence and Symbolic Computation, Proceedings of the 8th International Conference, AISC 2006, Beijing, China, September 20-22, 2006}, Springer (2006); \end{minipage}} \\ \\
\multicolumn{1}{l}{\begin{minipage}{0.1cm} \vskip-.6cm $\!\!$ \end{minipage}}
& \multicolumn{1}{c}{\begin{minipage}{14.3cm} \vskip-.4cm S. Liang, D.J. Jeffrey, and M.M. Maza {\it The Complete Root Classification of a Parametric
Polynomial on an Interval}, {\it ISSAC 2008, Proceedings of the 21st International Symposium on Symbolic and Algebraic Computation}, Linz/Hagenberg Austria (2010). \end{minipage}} \\ \\
\multicolumn{1}{l}{\begin{minipage}{0.1cm} \vskip-1.0cm $\!\!$[3] \end{minipage}}
& \multicolumn{1}{c}{\begin{minipage}{14.3cm} Emil M. Prodanov, {\it On the Determination of the Number of Positive and Negative Polynomial Zeros and Their Isolation}, Open Mathematics (De Gruyter) {\bf 18}, 1387--1412 (2020), arXiv: 1901.05960. \end{minipage}} \\ \\
\multicolumn{1}{l}{\begin{minipage}{0.1cm} \vskip-.5cm $\!\!$[4] \end{minipage}}
& \multicolumn{1}{c}{\begin{minipage}{14.3cm} Emil M. Prodanov, {\it Classification of the Roots of the Quartic Equation and their Pythagorean Tunes}, arXiv: 2008.07529 \end{minipage}} \\ \\
\multicolumn{1}{l}{\begin{minipage}{0.1cm} \vskip-.05cm $\!\!$[5] \end{minipage}}
& \multicolumn{1}{c}{\begin{minipage}{14.3cm} Emil M. Prodanov, {\it New Bounds on the Real Polynomial Roots}, arXiv:2008.11039. \end{minipage}} \\ \\
\multicolumn{1}{l}{\begin{minipage}{0.1cm} \vskip-.05cm $\!\!$[6] \end{minipage}}
& \multicolumn{1}{c}{\begin{minipage}{14.3cm} A.G. Kurosh, {\it Higher Algebra}, Mir Publishers Moscow (1972). \end{minipage}} \\ \\
\end{tabular}
\end{table}

\end{document}